\documentclass{kno_paper}

\title{Brownian Twins of Fractional Diffusion Models}
\author{%
  \knoauthor
}
\date{}

\usepackage[inline]{enumitem}
\setlist[enumerate,1]{label=(\roman*), ref=(\roman*)}

\newcommand{\GaussDensity}{\varphi}
\newcommand{\Smoothed}{p}
\newcommand{\Score}{s}
\newcommand{\vA}{v_{\mathrm{A}}}
\newcommand{\vB}{v_{\mathrm{B}}}
\newcommand{\Info}{\mathcal{F}}
\newcommand{\rev}[1]{\bar{#1}}
\newcommand{\Twin}{X^{\mathrm{B}}}
\newcommand{\LossTwin}{\mathbb{L}^{\mathrm{B}}}
\newcommand{\LossFrac}{\mathbb{L}}

\newcommand{\Id}{\IdMat_d}
\newcommand{\Nor}{\NormalDist}

\begin{document}

\maketitle

\begin{abstract}
    We study generative diffusion models in which the Brownian motion driving the noising dynamics is replaced by fractional Brownian motion or by a Markovian approximation of it.
    In the standard linear-Gaussian setting, we show that the memory of the noise acts only through the noise schedule.
    To every such model we associate a Brownian twin, a classical diffusion model with the same drift and an explicit variance schedule, and prove that everything the reverse-time process knows about the data is contained in one linear statistic which follows the reverse-time equation of the twin.
    Consequently, the denoising score matching losses of the two models coincide for every network, and access to the full path or to the augmented state does not reduce the loss.
    Moreover, we characterize when the reduction degenerates and describe where the reduction ends.

\end{abstract}

\section{Introduction}
\label{sec:introduction}

Diffusion models generate data by reversing a noising process~\cite{ho2020denoising,song2021score}.
In the standard construction, the noise is a Brownian motion, so that the noising process is Markovian and the reverse-time process needs only its current state and the score of the current marginal density.
Several recent works replace the Brownian motion by a noise with memory, namely by a fractional Brownian motion (fBM) or by a Markovian approximation of it through finitely many Ornstein--Uhlenbeck processes~\cite{daems2024variational,nobis2025generative.local,nobis2025fractional.local}.

In this paper, we investigate the effect of memory in the driving noise on these models.
In the standard linear-Gaussian setting, the answer is that it changes only the noise schedule.
More precisely, we consider forward processes $X$ and associate with each of them a classical diffusion model $\Twin$, which we call its \emph{Brownian twin}, where
\[
    X_t=c(t)\bigl(X_0+G_t\bigr)
    \qquad\text{and}\qquad
    \Twin_t\defined c(t)\bigl(X_0+W_{v(t)}\bigr),\qquad t\in[0,T].
\]
Here $X_0$ are the data, $c$ is a deterministic scale, $G$ is a Gaussian noise process independent of $X_0$, $W$ is a Brownian motion independent of $X_0$, and the variance schedule $v$ is computed from the covariance of $G$ alone.
This covers forward processes driven by fBM and the generative fractional diffusion models (GFDM) of \citet{nobis2025generative.local}.

The reduction rests on an elementary observation.
Whatever the reverse-time process has observed at time $t$, the noisy path on $[t,T]$ or an augmented state, the data enter these observations linearly and the noise is Gaussian.
We first show that such observations can be split into the best linear unbiased estimate $m_t$ of $X_0$ and a remainder which is independent of the data and of $m_t$, so that $m_t$ is a sufficient statistic and equals $X_0$ plus Gaussian noise of some variance $v(t)$.
Using the fact that the reverse-time process gains information as it moves towards time zero, we then show that the errors of these estimates are nested like the values of a Brownian motion.
From this, we are able to argue that \emph{the statistic $m_t$ follows the reverse-time equation of the Brownian twin}, and that it does so in the filtration of everything the reverse-time process has observed.
As a consequence, \emph{the denoising score matching losses of the fractional model and of its twin coincide for every network}, and no network with access to the full path or to the augmented state achieves a smaller loss.

We then make the reduction explicit.
For forward processes driven by fBM, the statistic is a linear filter of the noisy path, the schedule is given by the Molchan--Golosov kernel of the time-reversed noise, and we obtain the time reversal of the fractional path itself: The reversed noise is again an fBM plus a drift, and this drift is the classical score evaluated at the filter.
For GFDM, the statistic is the network input already used in~\cite{nobis2025generative.local}; here the reduction requires a nondegeneracy condition which can fail, and we characterize when it does.

Let us highlight what these results do and do not imply in practice.
A fractional model and its Brownian twin share the schedule, the regression problem solved in training, and the approximation made when the reverse-time dynamics are started from a Gaussian law at the final time.
They can still differ after discretization, since the two models integrate different equations, and in the optimization if the network is given a different input.
The Brownian model with the effective schedule is therefore the natural control experiment for any reported effect of memory.

\paragraph{Related Literature.}
Parts of the reduction are known for Markovian approximations.
For augmented linear diffusions in which the state given the data is Gaussian, \citet[Proposition~3.1]{chen2025tada} show that the posterior mean of the data depends on the augmented state only through the best linear unbiased estimator, so that the training objective is that of a conventional diffusion model; an equivalence of the dynamics is not stated there.
For GFDM, the fixed-time part of our statements is implicit in~\cite[Proposition~4.3 and~(17)]{nobis2025generative.local}: The score network takes the statistic as its input, and the scores of the auxiliary components are explicit functions of the learned one.
The bridge models of~\cite{nobis2025fractional.local} have the same structure.
In contrast to these works, our argument also covers genuine fBM, for which the reverse-time process observes a path instead of a finite-dimensional state.

Most closely related to our paper is \citet{han2026volterra}, who considers diffusion models driven by an exponential-sum approximation of a Volterra process, without drift and with isotropic noise.
There it is shown that the normalized residual has the path law of the data plus a time-changed Brownian motion, with a variance-preserving counterpart, and that pretrained denoisers can be reused at the effective noise level; see~\cite[Propositions~4, 6, and~7]{han2026volterra}.
In comparison, our argument applies to arbitrary linear-Gaussian observations and arbitrary data distributions, and it yields the reverse-time equation in the filtration of everything observed, which is stronger than an equality in law.
Moreover, for GFDM we include the drift, determine the effective schedule, and a criterion for degeneracy.

Filtering and estimation for signals observed in fBM noise go back to~\cite{lebreton1998filtering}.

The remainder of this paper is organized as follows.
In Section~\ref{sec:setting} we fix the setting and recall Brownian diffusion models with a given schedule.
Section~\ref{sec:reduction} contains the reduction in its abstract form.
In Sections~\ref{sec:fbm} and~\ref{sec:gfdm} we apply it to forward processes driven by fBM and to GFDM, respectively.
Finally, Section~\ref{sec:limits} describes the limits of the reduction, and Appendix~\ref{sec:auxiliary} and~\ref{sec:kernel} collect auxiliary results.

\paragraph{Notation.}
Throughout, we fix a complete probability space $(\Om,\Sig,\Prob)$, a time horizon $T>0$, and a dimension $d\in\N$, and all $\sigma$-fields and filtrations are augmented by the $\Prob$-null sets.
We denote by $\Abs{\argDot}$ the Euclidean norm on $\R^d$, by $\Id$ the identity matrix, and by $x^{(j)}$ the $j$th coordinate of $x\in\R^d$; scalars act on vectors coordinatewise.

\section{Setting}
\label{sec:setting}

In this section, we fix the data distribution and its Gaussian smoothings and recall Brownian diffusion models with a prescribed schedule, which serve as the reference class of the reduction.
We then introduce the common form of the fractional models studied in Sections~\ref{sec:fbm} and~\ref{sec:gfdm}.

\paragraph{Gaussian Smoothing and Scores.}
The \emph{data distribution} $\pi$ is a Borel probability measure on $\R^d$, and the data are a random vector $X_0\sim\pi$.
We assume neither a density nor finite moments, so that empirical distributions are included.
For $v>0$ we denote by $\GaussDensity_v$ the density of $\Nor(0,v\Id)$ and by $\xi$ a standard normal vector in $\R^d$ independent of $X_0$.
The random vector $X_0+\sqrt{v}\,\xi$ then has the density
\[
    \Smoothed_v(y)\defined\int_{\R^d}\GaussDensity_v(y-x)\,\pi(\de x),\qquad y\in\R^d,
\]
and we refer to $\Score_v\defined\nabla\log\Smoothed_v$ as the \emph{score at noise level $v$}.
The conditional law of $X_0$ given $X_0+\sqrt{v}\,\xi=y$ is the \emph{posterior}
\[
    \pi_{v,y}(\de x)\defined\frac{\GaussDensity_v(y-x)}{\Smoothed_v(y)}\,\pi(\de x).
\]
The following classical identity expresses the score through the posterior; see~\cite{robbins1956empirical,efron2011tweedie}.

\begin{lemma}[Tweedie's Formula]
    \label{lem:tweedie}
    The function $(v,y)\mapsto\Smoothed_v(y)$ is smooth and strictly positive on $(0,\infty)\times\R^d$.
    Moreover, for all $v>0$ and $y\in\R^d$, it holds that $\int\Abs{y-x}\,\pi_{v,y}(\de x)<\infty$ and
    \[
        \Score_v(y)=-\frac1v\int_{\R^d}(y-x)\,\pi_{v,y}(\de x).
    \]
\end{lemma}

\paragraph{Brownian Diffusion Models.}
The reference class consists of the classical diffusion models driven by Brownian motion~\cite{ho2020denoising,song2021score}.
We index these models by their schedule.

\begin{definition}
    \label{def:brownian-model}
    A \emph{schedule} is a pair $(c,v)$ of continuous functions $c\colon[0,T]\to(0,\infty)$ and $v\colon[0,T]\to[0,\infty)$ such that $c(0)=1$, $v(0)=0$, and $v$ is nondecreasing.
    The \emph{Brownian diffusion model with schedule $(c,v)$} is the process
    \[
        \Twin_t\defined c(t)\bigl(X_0+W_{v(t)}\bigr),\qquad t\in[0,T],
    \]
    where $W$ is a standard Brownian motion in $\R^d$ independent of $X_0$.
\end{definition}

Here $c$ scales the signal, and $v$ is the variance of the noise measured in units of the data.
If $c(t)=\exp(\int_0^t\mu(r)\de r)$ for some $\mu\in C([0,T])$ and $v$ is absolutely continuous, then $\Twin$ is a weak solution of
\begin{equation}
    \label{eq:twin-sde}
    \de\Twin_t=\mu(t)\Twin_t\de t+c(t)\sqrt{v'(t)}\de\widetilde{W}_t,\qquad t\in[0,T],
\end{equation}
for a Brownian motion $\widetilde{W}$.
The variance-exploding and the variance-preserving model of \citet{song2021score} correspond to $c\equiv1$ and to $v=c^{-2}-1$, respectively.
Given $X_0$, the state $\Twin_t$ has the law $\Nor(c(t)X_0,c(t)^2v(t)\Id)$.
For $v(t)>0$ it therefore has a density $p^{\mathrm{B}}_t$, and its score is given by
\begin{equation}
    \label{eq:twin-score}
    \nabla\log p^{\mathrm{B}}_t(x)=\frac{1}{c(t)}\,\Score_{v(t)}\bigl(x/c(t)\bigr),\qquad x\in\R^d.
\end{equation}

The score is learned by regression on the conditional score of $\Twin_t$ given $X_0$, which is explicit.
More precisely, let $v>0$ on $(0,T]$, let $\lambda\colon(0,T]\to[0,\infty)$ be a measurable weight, and let $\tau$ be a random time with values in $(0,T]$ which is independent of all other random variables.
The \emph{denoising score matching loss} of a measurable function $f\colon(0,T]\times\R^d\to\R^d$, which plays the role of the network, is
\begin{equation}
    \label{eq:dsm}
    \LossTwin(f)\defined\E\biggl[\lambda(\tau)\,\Bigl\lvert f(\tau,\Twin_\tau)+\frac{\Twin_\tau-c(\tau)X_0}{c(\tau)^2v(\tau)}\Bigr\rvert^2\biggr]\in[0,\infty].
\end{equation}
By the classical argument of \citet{vincent2011connection}, this loss is minimized by $f(t,\argDot)=\nabla\log p^{\mathrm{B}}_t$.
The minimal value may be infinite, which is why implementations bound $\tau$ away from zero or choose $\lambda$ accordingly; all identities below hold in $[0,\infty]$.
Noise prediction and data prediction are reparametrizations of the score, and all statements below about scores apply to them as well.

Generation runs in reverse time $s=T-t$, and we mark time-reversed objects with a bar.
In units of the data, the state of the Brownian model is $\Twin_t/c(t)=X_0+W_{v(t)}$.
If $v$ is absolutely continuous with $v'>0$ almost everywhere, then its time reversal $\rev{M}^{\mathrm{B}}_s\defined X_0+W_{v(T-s)}$ solves the classical reverse-time equation~\cite{anderson1982reverse,haussmann1986time}
\begin{equation}
    \label{eq:twin-reverse}
    \de\rev{M}^{\mathrm{B}}_s=v'(T-s)\,\Score_{v(T-s)}\bigl(\rev{M}^{\mathrm{B}}_s\bigr)\de s+\sqrt{v'(T-s)}\de\rev\beta_s,
    \qquad\rev{M}^{\mathrm{B}}_0\sim\Smoothed_{v(T)},
\end{equation}
for a Brownian motion $\rev\beta$ in the filtration generated by $\rev{M}^{\mathrm{B}}$.
In the present generality, this is a special case of Theorem~\ref{thm:twin} below.
A generative model simulates~\eqref{eq:twin-reverse}, or its rescaled version for $\Twin$, with a learned score in place of $\Score_v$.

\paragraph{Diffusion Models Driven by Gaussian Noise.}
The fractional models of Sections~\ref{sec:fbm} and~\ref{sec:gfdm} share a drift coefficient $\mu\in C([0,T])$ and a diffusion coefficient $g\in C^1([0,T])$ with $g>0$, and they differ in the noise.
We set
\[
    c(t)\defined\exp\Bigl(\int_0^t\mu(r)\de r\Bigr)
    \qquad\text{and}\qquad
    h\defined g/c\in C^1([0,T]).
\]
The forward equation and its solution then read
\begin{equation}
    \label{eq:forward}
    \de X_t=\mu(t)X_t\de t+g(t)\de N_t
    \quad\text{and}\quad
    X_t=c(t)\Bigl(X_0+\int_0^th(r)\de N_r\Bigr),\quad t\in[0,T],
\end{equation}
where the noise $N$ is a centered Gaussian process independent of $X_0$, namely a fractional Brownian motion in Section~\ref{sec:fbm} and a Markovian approximation of it in Section~\ref{sec:gfdm}.
The variance-exploding and variance-preserving choices of \citet{nobis2025generative.local} are of this form.
When the exact reverse-time process has reached time $t$, it has generated everything the model observes on $[t,T]$.
Since the data enter~\eqref{eq:forward} linearly and the noise is Gaussian, these observations are of the form studied in the next section.

\section{Reduction to a Brownian Model}
\label{sec:reduction}

In this section, we establish the reduction in an abstract form which covers both fractional models at once.
We first show that a family of linear-Gaussian observations carries the same information about the data as a single noisy copy of the data.
We then let the family grow, as it does for the reverse-time process on its way from time $T$ to time $0$, and show that these noisy copies evolve like a Brownian diffusion model run backwards.
Finally, we draw the consequences for training.

\subsection{Linear-Gaussian Observations}
\label{sec:observations}

Both fractional models observe the data through linear Gaussian channels.
We formulate this structure once.

\begin{definition}
    \label{def:observations}
    Let $I$ be an index set, and for every $i\in I$ let $a_i\in\R$ and let $G_i$ be a random vector in $\R^d$ such that
    \begin{enumerate}
        \item \label{obs:independent}
              the family $(G_i)_{i\in I}$ is independent of $X_0$, and
        \item \label{obs:gaussian}
              the coordinate families $(G^{(1)}_i)_{i\in I},\dots,(G^{(d)}_i)_{i\in I}$ are independent centered Gaussian families with the same covariance function.
    \end{enumerate}
    The random vectors $O_i\defined a_iX_0+G_i$, $i\in I$, are called \emph{linear-Gaussian observations} of $X_0$.
    For $J\subseteq I$ we set $\Info_J\defined\sigma(O_i:i\in J)$, and we say that $J$ \emph{sees the data} if $a_i\neq0$ for some $i\in J$.
\end{definition}

One may think of an index $i$ as a time $u\in[t,T]$, of $O_i$ as the observed state, of $a_i$ as the deterministic factor in front of $X_0$, and of $G_i$ as the noise.
For a finitely supported family $c=(c_i)_{i\in J}$ of real numbers we write $c\cdot O\defined\sum_ic_iO_i$ and $c\cdot G\defined\sum_ic_iG_i$.
The linear statistic $c\cdot O$ is called \emph{unbiased} if $\sum_ic_ia_i=1$.
In this case, $c\cdot O=X_0+c\cdot G$ is the data plus centered Gaussian noise, and we denote by
\[
    v_J\defined\inf\Bigl\{\tfrac1d\,\E\bigl[\Abs{c\cdot G}^2\bigr] : c\text{ finitely supported on }J\text{ with }\textstyle\sum_ic_ia_i=1\Bigr\}
\]
the smallest mean squared error per coordinate of an unbiased linear statistic.

The following result shows that all observations in $J$ together are worth exactly one noisy copy $m_J$ of the data with noise variance $v_J$.
In view of part~\ref{suff:blue}, we call $m_J$ the \emph{best linear unbiased estimator} of $X_0$ from the observations in $J$.
Everything else in the observations is noise which is independent of the data and of $m_J$.
This splitting is the key ingredient for all that follows.

\begin{proposition}[Sufficiency]
    \label{prop:sufficiency}
    Let $J\subseteq I$ see the data.
    Then there exists an $\Info_J$-measurable random vector $m_J$ with the following properties.
    \begin{enumerate}
        \item \label{suff:blue}
              If $(c_n)_{n\in\N}$ are unbiased with $\frac1d\E[\Abs{c_n\cdot G}^2]\to v_J$, then $\E[\Abs{c_n\cdot O-m_J}^2]\to0$.
        \item \label{suff:law}
              It holds that $m_J=X_0+e_J$, where $e_J\sim\Nor(0,v_J\Id)$ is independent of $X_0$.
        \item \label{suff:split}
              There is a $\sigma$-field $\mathcal{S}_J\subseteq\Info_J$, generated by Gaussian random variables, such that $\Info_J=\sigma(m_J)\vee\mathcal{S}_J$ and such that $X_0$, $e_J$, and $\mathcal{S}_J$ are independent.
        \item \label{suff:posterior}
              If $v_J>0$, then for every measurable $F\colon\R^d\times\R^d\to[0,\infty]$ it holds that
              \[
                  \E\bigl[F(X_0,m_J)\given\Info_J\bigr]=\int_{\R^d}F(x,m_J)\,\pi_{v_J,m_J}(\de x).
              \]
              In particular, the conditional law of $X_0$ given $\Info_J$ is the posterior $\pi_{v_J,m_J}$, and $\E[e_J\given\Info_J]=-v_J\,\Score_{v_J}(m_J)$.
              If $v_J=0$, then $X_0=m_J$ almost surely.
    \end{enumerate}
\end{proposition}

\begin{proof}
    By Definition~\ref{def:observations}~\ref{obs:gaussian}, the closed linear spans $\mathcal{H}^{(j)}$ of $\{G^{(j)}_i:i\in I\}$ in $\Lp^2(\Prob)$, $j=1,\dots,d$, are mutually orthogonal, and there are isometries $\iota_j\colon\mathcal{H}^{(1)}\to\mathcal{H}^{(j)}$ with $\iota_jG^{(1)}_i=G^{(j)}_i$.
    For $\zeta\in\mathcal{H}^{(1)}$ we write $\vec\zeta\defined(\iota_1\zeta,\dots,\iota_d\zeta)$, which is a random vector with independent $\Nor(0,\E[\zeta^2])$ coordinates.

    \emph{Step~1.} We construct $m_J$.
    Fix $i_0\in J$ with $a_{i_0}\neq0$ and consider the \emph{signal-free combinations}
    \[
        S_i\defined O_i-\frac{a_i}{a_{i_0}}\,O_{i_0}=G_i-\frac{a_i}{a_{i_0}}\,G_{i_0},\qquad i\in J,
    \]
    which are $\Info_J$-measurable and do not contain $X_0$.
    Let $\mathcal{H}_J\subseteq\mathcal{H}^{(1)}$ be the closed linear span of $\{S^{(1)}_i:i\in J\}$, let $P$ be the orthogonal projection onto $\mathcal{H}_J$, and set
    \[
        \zeta\defined\frac{G^{(1)}_{i_0}}{a_{i_0}},\qquad
        e\defined\zeta-P\zeta,\qquad
        e_J\defined\vec{e},\qquad
        m_J\defined\frac{O_{i_0}}{a_{i_0}}-\overrightarrow{P\zeta}.
    \]
    Since the maps $\iota_j$ are linear and satisfy $\iota_jG^{(1)}_{i_0}=G^{(j)}_{i_0}$, we have $\vec\zeta=G_{i_0}/a_{i_0}$ and $\vec\zeta-\overrightarrow{P\zeta}=\vec{e}$, and therefore
    \[
        m_J=X_0+\vec\zeta-\overrightarrow{P\zeta}=X_0+e_J.
    \]
    Each coordinate of $\overrightarrow{P\zeta}$ is an $\Lp^2$-limit of linear combinations of the variables $S^{(j)}_i$, $i\in J$, so $m_J$ is $\Info_J$-measurable.
    The coordinates of $e_J$ are independent and $\Nor(0,\E[e^2])$, and $e_J$ is a function of $(G_i)_{i\in I}$ and hence independent of $X_0$ by Definition~\ref{def:observations}~\ref{obs:independent}.

    \emph{Step~2.} We show that $\E[e^2]=v_J$ and prove~\ref{suff:blue}; together with Step~1, this also proves~\ref{suff:law}.
    Let $c$ be unbiased.
    Then $c\cdot G^{(1)}=\zeta+h$ with $h\defined\sum_ic_iS^{(1)}_i\in\mathcal{H}_J$, and since $e$ is orthogonal to $\mathcal{H}_J$, we obtain
    \begin{equation}
        \label{eq:pythagoras}
        \tfrac1d\,\E\bigl[\Abs{c\cdot G}^2\bigr]=\E\bigl[(c\cdot G^{(1)})^2\bigr]=\E[e^2]+\E\bigl[(P\zeta+h)^2\bigr]\geq\E[e^2].
    \end{equation}
    Conversely, every finite combination $h=\sum_id_iS^{(1)}_i$ arises in this way from an unbiased $c$, since $S_{i_0}=0$ leaves the coefficient $c_{i_0}$ free.
    Choosing $h$ close to $-P\zeta$ in $\Lp^2$ therefore shows that $\E[e^2]=v_J$.
    Now let $\frac1d\E[\Abs{c_n\cdot G}^2]\to v_J$, and let $h_n\in\mathcal{H}_J$ be associated with $c_n$ as above.
    Then $\E[(P\zeta+h_n)^2]\to0$ by~\eqref{eq:pythagoras}.
    Since $c_n\cdot O=X_0+c_n\cdot G$ and, by linearity of the $\iota_j$, $c_n\cdot G=\overrightarrow{\zeta+h_n}$, we obtain $c_n\cdot O-m_J=\overrightarrow{P\zeta+h_n}$, which tends to zero in $\Lp^2$ because the $\iota_j$ are isometries.
    In particular, $m_J$ does not depend on the choice of $i_0$.

    \emph{Step~3.} We prove~\ref{suff:split} with $\mathcal{S}_J\defined\sigma(S^{(j)}_i:i\in J,\ j=1,\dots,d)$, which is contained in $\Info_J$.
    Since $O_{i_0}/a_{i_0}=m_J+\overrightarrow{P\zeta}$ and $O_i=(a_i/a_{i_0})O_{i_0}+S_i$, every observation in $J$ is measurable with respect to $\sigma(m_J)\vee\mathcal{S}_J$, and hence $\Info_J=\sigma(m_J)\vee\mathcal{S}_J$.
    The $\sigma$-field $\mathcal{S}_J$ is generated by the subspace $\iota_1\mathcal{H}_J\oplus\dots\oplus\iota_d\mathcal{H}_J$.
    The coordinate $\iota_je$ of $e_J$ is orthogonal to $\iota_j\mathcal{H}_J$, because $\iota_j$ is an isometry, and to $\iota_k\mathcal{H}_J$ for $k\neq j$, because different coordinates are independent.
    Lemma~\ref{lem:gaussian-hilbert}~\ref{gh:independent} thus shows that $e_J$ and $\mathcal{S}_J$ are independent.
    Both are functions of $(G_i)_{i\in I}$ and hence independent of $X_0$, so $X_0$, $e_J$, and $\mathcal{S}_J$ are independent.

    \emph{Step~4.} We prove~\ref{suff:posterior}.
    By~\ref{suff:split}, the pair $(X_0,m_J)$ is independent of $\mathcal{S}_J$, and the sets $D_1\cap D_2$ with $D_1\in\sigma(m_J)$ and $D_2\in\mathcal{S}_J$ form a generator of $\Info_J$ which is stable under intersections.
    It follows that $\E[F(X_0,m_J)\given\Info_J]=\E[F(X_0,m_J)\given m_J]$, first for bounded $F$ and then, by monotone convergence, for all measurable $F\geq0$.
    If $v_J>0$, then $(X_0,m_J)$ has the law $\GaussDensity_{v_J}(m-x)\,\pi(\de x)\de m$, so the conditional law of $X_0$ given $m_J$ is $\pi_{v_J,m_J}$.
    Applying this to the positive and negative parts of the coordinates of $m_J-X_0$ and using Lemma~\ref{lem:tweedie}, we obtain
    \[
        \E\bigl[e_J\given\Info_J\bigr]=\int_{\R^d}(m_J-x)\,\pi_{v_J,m_J}(\de x)=-v_J\,\Score_{v_J}(m_J).
    \]
    If $v_J=0$, then $e_J=0$ almost surely, and the proof is complete.
\end{proof}

\subsection{Nested Observations and the Brownian Twin}
\label{sec:nested}

We now let the set of observations grow.
Let $(I_t)_{t\in(0,T]}$ be subsets of $I$ such that $I_{t'}\supseteq I_t$ for $t'\leq t$ and such that $I_T$ sees the data, and write
\[
    \Info_t\defined\Info_{I_t},\qquad
    m_t\defined m_{I_t},\qquad
    e_t\defined e_{I_t},\qquad
    v(t)\defined v_{I_t},\qquad t\in(0,T].
    \]
    In reverse time, we set
    \[
    \rev{M}_s\defined m_{T-s},\quad
    \rev\Info_s\defined\Info_{T-s}\quad
    \text{for }s\in[0,T),\qquad
    \rev{M}_T\defined X_0,\quad
    \rev\Info_T\defined\bigvee_{s<T}\rev\Info_s,
\]
so that $\rev\Info$ is a filtration.
Under the hypothesis of part~\ref{twin:sde} below, $v(t)\to0$ as $t\downarrow0$, and the value $\rev{M}_T=X_0$ is the limit of $\rev{M}_s$ as $s\uparrow T$; equivalently, $X_0=m_{I_0}$ almost surely for $I_0\defined\bigcup_{t>0}I_t$ by Proposition~\ref{prop:sufficiency}, since $v_{I_0}\leq v(t)$ for every $t>0$.
As the reverse-time process observes more, its estimate of the data improves, and the next theorem, which is the main result of this section, shows that the errors of the estimates are nested like the values of a Brownian motion.
A notable feature is that the reverse-time equation in part~\ref{twin:sde} holds in the filtration of everything observed and not only in the filtration generated by the estimates: Knowing the whole observed path does not change the drift, which is the classical score of the data distribution evaluated at the estimate.

\begin{theorem}[Brownian Twin]
    \label{thm:twin}
    In the situation described above, the following hold.
    \begin{enumerate}
        \item \label{twin:law}
              It holds that $\Cov[e_t,e_{t'}]=v(t\wedge t')\Id$ for all $t,t'\in(0,T]$, and $v$ is nondecreasing.
              Consequently, if $W$ is a standard Brownian motion in $\R^d$ independent of $X_0$, then $(X_0,(m_t)_{t\in(0,T]})$ and $(X_0,(X_0+W_{v(t)})_{t\in(0,T]})$ have the same finite-dimensional distributions.
        \item \label{twin:sde}
              Suppose that $v$ extends to an absolutely continuous function on $[0,T]$ with $v(0)=0$ and that $v>0$ on $(0,T]$.
              Then $m$ can be chosen to be continuous with $m_t\to X_0$ as $t\downarrow0$.
              Moreover, the process $\rev{M}$ is $\rev\Info$-adapted, it holds that $\E\int_0^Tv'(T-r)\Abs{\Score_{v(T-r)}(\rev{M}_r)}\de r<\infty$, and
              \[
                  \rev{N}_s\defined\rev{M}_s-\rev{M}_0-\int_0^sv'(T-r)\,\Score_{v(T-r)}\bigl(\rev{M}_r\bigr)\de r,\qquad s\in[0,T],
              \]
              is a continuous $\rev\Info$-martingale with $\langle\rev{N}^{(i)},\rev{N}^{(j)}\rangle_s=\delta_{ij}(v(T)-v(T-s))$.
              If, in addition, $v'>0$ almost everywhere, then
              \[
                  \rev\beta_s\defined\int_0^sv'(T-r)^{-1/2}\de\rev{N}_r,\qquad s\in[0,T],
              \]
              is an $\rev\Info$-Brownian motion, and $\rev{M}$ solves the reverse-time equation~\eqref{eq:twin-reverse} of the Brownian model with schedule $(1,v)$ in the filtration $\rev\Info$.
    \end{enumerate}
\end{theorem}

\begin{proof}
    We use the same $i_0\in I_T$ with $a_{i_0}\neq0$ in the construction of Proposition~\ref{prop:sufficiency} for every $t$.
    We write $\mathcal{H}_t\defined\mathcal{H}_{I_t}$, $P_t$ for the projection onto $\mathcal{H}_t$, and $e^{(1)}_t=\zeta-P_t\zeta$ for the first coordinate of $e_t$, and we note that $\mathcal{H}_{t'}\supseteq\mathcal{H}_t$ for $t'\leq t$.

    \emph{Step~1.} We prove~\ref{twin:law}.
    For $t'\leq t$, the difference $e^{(1)}_t-e^{(1)}_{t'}=P_{t'}\zeta-P_t\zeta$ lies in $\mathcal{H}_{t'}$ and is therefore orthogonal to $e^{(1)}_{t'}$.
    Hence $\E[e^{(1)}_te^{(1)}_{t'}]=\E[(e^{(1)}_{t'})^2]=v(t')$, which gives the covariance since the coordinates are independent copies, and $v(t')^2\leq v(t)v(t')$ by the Cauchy--Schwarz inequality.
    The processes $(e_t)_{t\in(0,T]}$ and $(W_{v(t)})_{t\in(0,T]}$ are thus centered Gaussian with the same covariance and independent of $X_0$, which implies the claim on the finite-dimensional distributions.

    \emph{Step~2.} We construct the continuous modification of $m$.
    Let $\mathbb{D}\subseteq(0,T]$ be countable and dense.
    Since $t\mapsto W_{v(t)}$ is continuous on $[0,T]$ and vanishes at zero, its restriction to $\mathbb{D}$ is almost surely uniformly continuous and tends to zero at zero.
    This is an event determined by countably many coordinates, so by~\ref{twin:law} the same holds for $(e_t)_{t\in\mathbb{D}}$.
    The continuous extension $\tilde{e}$ of $(e_t)_{t\in\mathbb{D}}$ to $[0,T]$ satisfies $\tilde{e}_t=e_t$ almost surely for every $t\in(0,T]$, because $\E\Abs{e_r-e_t}^2=d\Abs{v(r)-v(t)}\to0$ as $r\to t$.
    Hence $X_0+\tilde{e}$ is a continuous modification of $m$, and it is adapted to $(\Info_t)_{t\in(0,T]}$ since these $\sigma$-fields are complete.
        We write $e$ and $m$ for the modifications and set $\rev{e}_s\defined e_{T-s}$, so that $\rev{M}=X_0+\rev{e}$ is continuous on $[0,T]$ and $\rev\Info$-adapted on $[0,T)$.
    At $s=T$, the value $\rev{M}_T=X_0=\lim_{s\uparrow T}\rev{M}_s$ is $\rev\Info_T$-measurable since $\rev\Info_T$ contains all $\rev\Info_s$, $s<T$, and is complete.

    \emph{Step~3.} We show that $\rev{e}$ is a Gaussian bridge in the filtration $\FiltAlt_s\defined\rev\Info_s\vee\sigma(X_0)$.
    Let $\rho(s)\defined v(T-s)$, let $0\leq s\leq s'<T$, and put $t\defined T-s$ and $t'\defined T-s'$.
    By Proposition~\ref{prop:sufficiency}~\ref{suff:split}, we have $\FiltAlt_s=\sigma(X_0)\vee\sigma(e_t)\vee\mathcal{S}_{I_t}$.
    As in Step~3 of the proof of Proposition~\ref{prop:sufficiency}, the $\sigma$-field $\sigma(e_t)\vee\mathcal{S}_{I_t}$ is generated by the closed subspace $\mathcal{K}_t$ of $\Lp^2(\Prob)$ spanned by the coordinates $\iota_1e^{(1)}_t,\dots,\iota_de^{(1)}_t$ of $e_t$ and by the spaces $\iota_1\mathcal{H}_t,\dots,\iota_d\mathcal{H}_t$.
    Since $e^{(1)}_t$ is orthogonal to $\mathcal{H}_t$ and different coordinates are independent, $\mathcal{K}_t$ is the orthogonal sum of the spaces $\iota_k\mathcal{H}_t$ and of the one-dimensional spaces spanned by $\iota_ke^{(1)}_t$, $k=1,\dots,d$.
    We fix a coordinate $j$ and compute the orthogonal projection of $\iota_je^{(1)}_{t'}$, the $j$th coordinate of $e_{t'}$, onto each of these spaces.
    The variable $\iota_je^{(1)}_{t'}$ is orthogonal to $\iota_j\mathcal{H}_t$, because $e^{(1)}_{t'}$ is orthogonal to $\mathcal{H}_{t'}\supseteq\mathcal{H}_t$, and it is orthogonal to all spaces with index $k\neq j$, because different coordinates are independent.
    By~\ref{twin:law}, the projection of $\iota_je^{(1)}_{t'}$ onto the span of $\iota_je^{(1)}_t$ is
    \[
        \frac{\E\bigl[e^{(1)}_{t'}e^{(1)}_t\bigr]}{\E\bigl[(e^{(1)}_t)^2\bigr]}\,\iota_je^{(1)}_t=\frac{v(t')}{v(t)}\,\iota_je^{(1)}_t.
    \]
    Hence the projection of $\iota_je^{(1)}_{t'}$ onto $\mathcal{K}_t$ is $(v(t')/v(t))\,\iota_je^{(1)}_t$.
    Since $X_0$ is independent of all Gaussian variables involved, Lemma~\ref{lem:gaussian-hilbert}~\ref{gh:projection} yields
    \[
        \E\bigl[\rev{e}_{s'}\given\FiltAlt_s\bigr]=\frac{v(t')}{v(t)}\,e_t=\frac{\rho(s')}{\rho(s)}\,\rev{e}_s.
    \]
    Moreover, by~\ref{twin:law}, the coordinates of $\rev{e}$ are independent centered Gaussian processes with covariance $\rho(s\vee s')$.
    Lemma~\ref{lem:bridge} in the appendix therefore applies to $\rev{e}$ and shows that
    \[
        A_s\defined\rev{e}_s-\rev{e}_0-\int_0^sb_r\de r
        \qquad\text{with}\qquad
        b_r\defined\frac{\rho'(r)}{\rho(r)}\,\rev{e}_r
    \]
    is a continuous $\FiltAlt$-martingale with $\langle A^{(i)},A^{(j)}\rangle_s=\delta_{ij}(\rho(0)-\rho(s))$ and that $\E\int_0^T\Abs{b_r}\de r<\infty$.

    \emph{Step~4.} We project the drift onto the observed information.
    We have $\rev{M}_s-\rev{M}_0=\int_0^sb_r\de r+A_s$ with an $\rev\Info$-adapted left-hand side, and Proposition~\ref{prop:sufficiency}~\ref{suff:posterior} with $J=I_{T-r}$ gives
    \[
        \E\bigl[b_r\given\rev\Info_r\bigr]
        =\frac{\rho'(r)}{\rho(r)}\,\E\bigl[e_{T-r}\given\Info_{T-r}\bigr]
        =-\rho'(r)\,\Score_{\rho(r)}(\rev{M}_r)
        =v'(T-r)\,\Score_{v(T-r)}(\rev{M}_r)
    \]
    for almost every $r\in[0,T)$.
    The right-hand side is $\rev\Info$-progressively measurable, because $(v,y)\mapsto\Score_v(y)$ is continuous by Lemma~\ref{lem:tweedie} and $\rev{M}$ is continuous and adapted, and it is integrable by Jensen's inequality.
    Lemma~\ref{lem:innovation} now shows that $\rev{N}$ is an $\rev\Info$-martingale.
    Since $\rev{N}-A$ is continuous and of finite variation and quadratic covariations do not depend on the filtration, $\rev{N}$ has the same brackets as $A$, and $\rho(0)-\rho(s)=v(T)-v(T-s)$.

    \emph{Step~5.} We identify the Brownian motion.
    If $v'>0$ almost everywhere, then $\rev\beta$ is a continuous $\rev\Info$-local martingale with $\langle\rev\beta^{(i)},\rev\beta^{(j)}\rangle_s=\delta_{ij}s$ and hence an $\rev\Info$-Brownian motion by L\'evy's characterization.
    Since $\rev{N}_s=\int_0^s\sqrt{v'(T-r)}\de\rev\beta_r$ and $\rev{M}_0=X_0+e_T$ has the density $\Smoothed_{v(T)}$, this is~\eqref{eq:twin-reverse}, and the proof is complete.
\end{proof}

\subsection{Consequences for Training}
\label{sec:training}

Let $c\colon[0,T]\to(0,\infty)$ be continuous with $c(0)=1$, and suppose that the hypotheses of Theorem~\ref{thm:twin}~\ref{twin:sde} hold, so that $(c,v)$ is a schedule.
We call the Brownian diffusion model $\Twin$ with schedule $(c,v)$ the \emph{Brownian twin} of the observation model.
Moreover, we refer to
\[
    \tilde{x}_t\defined c(t)m_t
    \qquad\text{and}\qquad
    \Lambda_t\defined-\frac{\tilde{x}_t-c(t)X_0}{c(t)^2v(t)},\qquad t\in(0,T],
\]
as the \emph{network input} and the \emph{regression target}, respectively.
In analogy to~\eqref{eq:dsm}, the denoising score matching loss of the observation model is
\begin{equation}
    \label{eq:loss-frac}
    \LossFrac(f)\defined\E\bigl[\lambda(\tau)\,\Abs{f(\tau,\tilde{x}_\tau)-\Lambda_\tau}^2\bigr]\in[0,\infty]
\end{equation}
for measurable $f\colon(0,T]\times\R^d\to\R^d$.
The following corollary shows that a network trained on the observation model and a network trained on its Brownian twin solve the same regression problem, and that no additional input improves it.

\begin{corollary}[Score Matching]
    \label{cor:score-matching}
    In the situation described above, the following hold.
    \begin{enumerate}
        \item \label{sm:loss}
              The pairs $(X_0,\tilde{x})$ and $(X_0,\Twin)$ have the same law.
              Consequently, $\E[\ell(\tau,X_0,\tilde{x}_\tau)]=\E[\ell(\tau,X_0,\Twin_\tau)]$ for every measurable $\ell\geq0$, and in particular $\LossFrac(f)=\LossTwin(f)$ for every $f$, every weight $\lambda$, and every law of $\tau$.
        \item \label{sm:minimizer}
              Both losses are minimized by the score of the twin,
              \[
                  f^*(t,x)\defined\nabla\log p^{\mathrm{B}}_t(x)=\frac{1}{c(t)}\,\Score_{v(t)}\bigl(x/c(t)\bigr),
              \]
              and if the minimum is finite, then every minimizer agrees with $f^*$ almost everywhere on $\{\lambda>0\}$ with respect to the law of $(\tau,\Twin_\tau)$.
        \item \label{sm:no-gain}
              For every $t\in(0,T]$ and every $\Info_t$-measurable random vector $\Phi$, it holds that
              \[
                  \E\Big[\Abs{\Phi-\Lambda_t}^2\Big]\geq\E\Big[\Abs{f^*(t,\tilde{x}_t)-\Lambda_t}^2\Big],
              \]
              with equality if and only if $\Phi=f^*(t,\tilde{x}_t)$ almost surely.
    \end{enumerate}
\end{corollary}

\begin{proof}
    The first claim in~\ref{sm:loss} is Theorem~\ref{thm:twin}~\ref{twin:law}, multiplied by $c(t)$, and using that the law of a process is determined by its finite-dimensional distributions.
    For the second claim, we note that $\tilde{x}$ and $\Twin$ are continuous on $(0,T]$ and hence jointly measurable, and that $\tau$ is independent of all other random variables.
    Writing $\vartheta$ for the law of $\tau$, Fubini's theorem therefore gives
    \[
        \E\bigl[\ell(\tau,X_0,\tilde{x}_\tau)\bigr]
        =\int\E\bigl[\ell(t,X_0,\tilde{x}_t)\bigr]\,\vartheta(\de t)
        =\int\E\bigl[\ell(t,X_0,\Twin_t)\bigr]\,\vartheta(\de t)
        =\E\bigl[\ell(\tau,X_0,\Twin_\tau)\bigr].
    \]
    To prove~\ref{sm:no-gain}, we observe that $\Lambda_t=-e_t/(c(t)v(t))$ is square integrable, and that Proposition~\ref{prop:sufficiency}~\ref{suff:posterior} with $J=I_t$ and~\eqref{eq:twin-score} show that
    \[
        \E\bigl[\Lambda_t\given\Info_t\bigr]=-\frac{\E[e_t\given\Info_t]}{c(t)v(t)}=\frac{\Score_{v(t)}(m_t)}{c(t)}=f^*(t,\tilde{x}_t).
    \]
    Hence, for every square integrable $\Info_t$-measurable $\Phi$,
    \[
        \E\Big[\Abs{\Phi-\Lambda_t}^2\Big]=\E\Big[\Abs{\Phi-f^*(t,\tilde{x}_t)}^2\Big]+\E\Big[\Abs{f^*(t,\tilde{x}_t)-\Lambda_t}^2\Big],
    \]
    while the left-hand side is infinite if $\Phi$ is not square integrable.
    This is~\ref{sm:no-gain}.
    Finally, applying this identity to $\Phi=f(t,\tilde{x}_t)$ and integrating against $\lambda(t)\,\vartheta(\de t)$ yields
    \[
        \LossFrac(f)=\E\bigl[\lambda(\tau)\,\Abs{f(\tau,\tilde{x}_\tau)-f^*(\tau,\tilde{x}_\tau)}^2\bigr]+\LossFrac(f^*),
    \]
    which proves~\ref{sm:minimizer} for $\LossFrac$ and, by~\ref{sm:loss}, for $\LossTwin$.
\end{proof}

It follows, in particular, that the continuous processes $\rev{M}$ and $(\Twin_{T-s}/c(T-s))_{s\in[0,T]}$ have the same law, and both end in $X_0\sim\pi$.
With the exact score, the observation model and its Brownian twin thus generate the data distribution along the same reverse-time dynamics.
It remains to identify the statistic $m_t$ and the schedule $v$ for the two fractional models, which is the content of the following two sections.

\section{Forward Processes Driven by Fractional Brownian Motion}
\label{sec:fbm}

In this section, the noise in~\eqref{eq:forward} is a fractional Brownian motion.
We show that the reverse-time process, although it observes the whole path up to a given point in time, sees the data only through one linear statistic, and we compute the schedule of its Brownian twin explicitly.
The key ingredient is the Molchan--Golosov representation of the reversed noise, which turns the observed path into the filtration of a Brownian motion; an initial enlargement of this filtration then yields the reverse-time dynamics.

Let $H\in(0,1)$.
Following \citet{mandelbrot1968fractional}, a \emph{fractional Brownian motion} (fBM) with Hurst index $H$ is a centered real Gaussian process $B=(B_t)_{t\geq0}$ with continuous paths and covariance
\[
    \E[B_sB_t]=R_H(s,t)\defined\tfrac12\bigl(s^{2H}+t^{2H}-\Abs{t-s}^{2H}\bigr),\qquad s,t\geq0.
\]
An fBM in $\R^d$ consists of independent fBM coordinates with the same Hurst index.
For $H=\frac12$ this is Brownian motion, while for $H\neq\frac12$ the increments are correlated, positively for $H>\frac12$ and negatively for $H<\frac12$.
Its increments are stationary with $\E[(B_a-B_b)^2]=\Abs{a-b}^{2H}$, so polarization shows that the reversed increments $\rev{B}_s\defined B_{T-s}-B_T$, $s\in[0,T]$, form again an fBM with the same Hurst index.

We use the Volterra representation of fBM with respect to a Brownian motion which generates the same filtration.
It can be found in~\cite{decreusefond1999stochastic,molchan1969gaussian,norros1999elementary} and~\cite{nualart2006malliavin}, and the equality of the filtrations is stated explicitly in~\cite{picard2011representation}.

\begin{lemma}[Molchan--Golosov Representation]
    \label{lem:molchan-golosov}
    There is a measurable kernel $K_H\colon\{(t,r):0<r<t\}\to(0,\infty)$ with
    \[
        \int_0^{t\wedge u}K_H(t,r)K_H(u,r)\de r=R_H(t,u),\qquad t,u>0,
    \]
    in particular $\int_0^tK_H(t,r)^2\de r=t^{2H}$, such that the following holds.
    For every fBM $B$ with Hurst index $H$ on $[0,T]$ there is a standard Brownian motion $W$ on $[0,T]$, defined on the same probability space, with
    \[
        B_t=\int_0^tK_H(t,r)\de W_r
        \qquad\text{and}\qquad
        \Filt^B_t=\Filt^W_t,\qquad t\in[0,T],
    \]
    where $\Filt^B$ and $\Filt^W$ denote the augmented natural filtrations.
    For $H=\frac12$, we have $K_H\equiv1$ and $W=B$.
\end{lemma}

The kernel is given in Appendix~\ref{sec:kernel}, where Lemma~\ref{lem:kernel} collects the properties we need.
Integrals of deterministic functions with respect to fBM are defined by integration by parts, that is, for $f\in C^1([0,T])$ and an fBM $B$ we set
\begin{equation}
    \label{eq:wiener-integral}
    \int_0^tf(r)\de B_r\defined f(t)B_t-\int_0^tf'(r)B_r\de r,
    \qquad t\in[0,T].
\end{equation}
For $H=\frac12$ this is the It\^{o} integral, for $H>\frac12$ the pathwise Riemann--Stieltjes integral, and for all $H$ it agrees with the Wiener integral of $f\Ind{[0,t]}$ with respect to $B$; see~\cite{nualart2006malliavin}.

Let $B$ be an fBM in $\R^d$ with Hurst index $H\in(0,1)$, independent of $X_0$.
With the coefficients $\mu$, $g$, $c$, and $h=g/c$ of Section~\ref{sec:setting}, the forward process is~\eqref{eq:forward} with $N=B$, that is,
\[
    X_t=c(t)\bigl(X_0+G_t\bigr)
    \qquad\text{with}\qquad
    G_t\defined\int_0^th(r)\de B_r,\qquad t\in[0,T].
\]
When the exact reverse-time process has reached time $t$, it has generated the path of $X$ on $[t,T]$, so its information is
\[
    \Info^{\mathrm{A}}_t\defined\sigma\bigl(X_u:u\in[t,T]\bigr),\qquad t\in[0,T].
\]
For $H\neq\frac12$, the Markov property is not available, so a priori the reverse-time process depends on the whole path.

In order to state the main result of this section, we need to introduce some notation.
First, we denote by $\rev{W}$ the Brownian motion of the Molchan--Golosov representation of the reversed fBM $\rev{B}$, so that
\[
    \rev{B}_s=\int_0^sK_H(s,r)\de\rev{W}_r
    \qquad\text{and}\qquad
    \Filt^{\rev{B}}_s=\Filt^{\rev{W}}_s,\qquad s\in[0,T].
\]
Second, the \emph{reversed noise} is given by
\begin{equation}
    \label{eq:reversed-noise}
    \rev{G}_s\defined G_{T-s}-G_T=\int_0^sh(T-r)\de\rev{B}_r,
    \qquad s\in[0,T],
\end{equation}
where the integral is understood in the sense of~\eqref{eq:wiener-integral}; we verify this identity in the proof below.
Third, we define the \emph{terminal kernel}
\begin{equation}
    \label{eq:k-h}
    k_h(r)\defined h(0)K_H(T,r)+\int_r^Th'(T-u)K_H(u,r)\de u,
    \qquad r\in(0,T).
\end{equation}
We will see that $k_h\in\Lp^2(0,T)$ and that $\rev{G}_T=-G_T=\int_0^Tk_h\de\rev{W}$.
With this notation in place, the main result of this section can be stated as follows.

\begin{theorem}[Effective Schedule of Fractional Noise]
    \label{thm:fbm}
    The following hold.
    \begin{enumerate}
        \item \label{fbm:statistic}
              For $t\in(0,T]$, the random vector
              \[
                  m_t\defined\frac{X_T}{c(T)}-\E\bigl[G_T\given G_u-G_T,\ u\in[t,T]\bigr]
              \]
              is the best linear unbiased estimator of $X_0$ from the path of $X$ on $[t,T]$.
              That is, in the notation of Theorem~\ref{thm:twin} with $I_t\defined[t,T]$ and the observations $X_u/c(u)$, $u\in[0,T]$, we have $m_t=m_{I_t}$ and $v(t)=\vA(t)$, where
              \[
                  \vA(t)\defined\int_{T-t}^Tk_h(r)^2\de r>0.
              \]
              Moreover, $\vA$ with $\vA(0)\defined0$ is absolutely continuous on $[0,T]$ with
              \[
                  \vA'(t)=k_h(T-t)^2,
                  \qquad
                  \vA(T)=\Var\bigl[G^{(1)}_T\bigr],
                  \qquad\text{and}\qquad
                  \vA(t)\leq\Var\bigl[G^{(1)}_t\bigr],
              \]
              and if $h'\geq0$, then $k_h>0$ on $(0,T)$, so that $\vA'>0$ almost everywhere.
        \item \label{fbm:statistic-sde}
              In reverse time, the statistic is a Wiener integral with respect to $\rev{W}$, namely
              \[
                  \rev{M}_s\defined m_{T-s}=\frac{X_T}{c(T)}+\int_0^sk_h(r)\de\rev{W}_r,
                  \qquad s\in[0,T).
              \]
        \item \label{fbm:path-sde}
              Let $\rev{v}(s)\defined\vA(T-s)$.
              The process
              \[
                  \widehat{W}_s\defined\rev{W}_s-\int_0^sk_h(r)\,\Score_{\rev{v}(r)}\bigl(\rev{M}_r\bigr)\de r,
                  \qquad s\in[0,T),
              \]
              is a Brownian motion in the filtration $(\Info^{\mathrm{A}}_{T-s})_{s\in[0,T)}$ and extends continuously to a Brownian motion on $[0,T]$.
              Consequently, for $s\in[0,T)$, the reversed path $\rev{X}_s\defined X_{T-s}/c(T-s)$ satisfies
              \begin{align}
                  \rev{X}_s & =\rev{X}_0+\int_0^sh(T-r)\de\rev{B}_r,                                                               \\
                  \rev{B}_s & =\int_0^sK_H(s,r)\Bigl(\de\widehat{W}_r+k_h(r)\,\Score_{\rev{v}(r)}\bigl(\rev{M}_r\bigr)\de r\Bigr),
              \end{align}
              and the martingale part of $\rev{M}$ is $\int_0^{\argDot}k_h\de\widehat{W}$.
    \end{enumerate}
    In particular, the hypotheses of Theorem~\ref{thm:twin}~\ref{twin:sde} hold, so that Proposition~\ref{prop:sufficiency}, Theorem~\ref{thm:twin}, and Corollary~\ref{cor:score-matching} apply with $v=\vA$ and $\Info_t=\Info^{\mathrm{A}}_t$.
\end{theorem}

\begin{proof}
    We first express the information $\Info^{\mathrm{A}}_t$ and the terminal noise $G_T$ through the Brownian motion $\rev{W}$.
    Parts~\ref{fbm:statistic} and~\ref{fbm:statistic-sde} then follow from Proposition~\ref{prop:sufficiency}, and part~\ref{fbm:path-sde} follows from an initial enlargement of the filtration of $\rev{W}$.

    \emph{Step~1.}
    We verify~\eqref{eq:reversed-noise} and identify the information $\Info^{\mathrm{A}}_t$.
    By~\eqref{eq:wiener-integral} and $B_r=\rev{B}_{T-r}+B_T$, we have
    \begin{align}
        \rev{G}_s = G_{T-s}-G_T
         & =h(T-s)B_{T-s}-h(T)B_T+\int_{T-s}^Th'(r)B_r\de r                \\
         & =h(T-s)\rev{B}_s+\int_0^sh'(T-q)\,\rev{B}_q\de q                \\
         & \hspace{3cm}+B_T\Bigl(h(T-s)-h(T)+\int_{T-s}^Th'(r)\de r\Bigr),
    \end{align}
    and the last bracket vanishes.
    With $\rev{h}(q)\defined h(T-q)$, so that $\rev{h}'(q)=-h'(T-q)$, we conclude
    \[
        \rev{G}_s=G_{T-s}-G_T=h(T-s)\rev{B}_s+\int_0^sh'(T-q)\,\rev{B}_q\de q=\rev{h}(s)\rev{B}_s-\int_0^s\rev{h}'(q)\rev{B}_q\de q,
    \]
    which is~\eqref{eq:reversed-noise}.
    Conversely, the reversed noise determines the reversed fBM.
    Since $h$ is positive and continuously differentiable, dividing the identity above by $\rev{h}(s)$ shows that every path of $\rev{B}$ solves the linear Volterra equation of the second kind
    \[
        \rev{B}_s=\frac{\rev{G}_s}{\rev{h}(s)}+\int_0^s\frac{\rev{h}'(q)}{\rev{h}(s)}\,\rev{B}_q\de q,\qquad s\in[0,T],
    \]
    whose kernel is continuous on $\{0\leq q\leq s\leq T\}$.
    This equation has a unique continuous solution, the uniform limit of the Picard iterates started at $\rev{G}/\rev{h}$, and uniqueness follows from Gronwall's inequality.
    The value of each iterate at time $s$ is a measurable functional of $(\rev{G}_q)_{q\leq s}$, so $\rev{B}_s$ is $\Filt^{\rev{G}}_s$-measurable for every $s\in[0,T]$.
    As $\rev{G}_s$ is in turn a functional of $(\rev{B}_q)_{q\leq s}$ by the identity above, we obtain $\Filt^{\rev{G}}_s=\Filt^{\rev{B}}_s$, and $\Filt^{\rev{B}}_s=\Filt^{\rev{W}}_s$ from Lemma~\ref{lem:molchan-golosov}.
    Since $X_u/c(u)=X_T/c(T)+(G_u-G_T)$ and $\sigma(G_u-G_T:u\in[t,T])=\Filt^{\rev{G}}_{T-t}$, we obtain
    \begin{equation}
        \label{eq:info-a}
        \Info^{\mathrm{A}}_t=\sigma\bigl(X_T/c(T)\bigr)\vee\Filt^{\rev{W}}_{T-t},
        \qquad t\in[0,T].
    \end{equation}

    \emph{Step~2.}
    We show that $k_h\in\Lp^2(0,T)$ with $\rev{G}_T=\int_0^Tk_h\de\rev{W}$, and that $k_h$ is positive near $T$.
    The map $u\mapsto\rev{B}_u=\int_0^uK_H(u,r)\de\rev{W}_r$ is continuous from $[0,T]$ into $\Lp^2(\Prob)$.
    By Lemma~\ref{lem:molchan-golosov}, we have $\Norm{K_H(u,\argDot)}_{\Lp^2}=u^H$, so that $\int_0^T\Abs{h'(T-u)}\,\Norm{K_H(u,\argDot)}_{\Lp^2}\de u<\infty$.
    Since the Wiener integral is an isometry from $\Lp^2(0,T)$ into $\Lp^2(\Prob)$, it commutes with the corresponding Bochner integral.
    Using~\eqref{eq:reversed-noise} for $s=T$ and~\eqref{eq:wiener-integral}, we therefore obtain
    \begin{align}
        \rev{G}_T
         & =h(0)\rev{B}_T+\int_0^Th'(T-u)\,\rev{B}_u\de u                             \\
         & =\int_0^T\Bigl(h(0)K_H(T,r)+\int_r^Th'(T-u)K_H(u,r)\de u\Bigr)\de\rev{W}_r
        =\int_0^Tk_h(r)\de\rev{W}_r,
    \end{align}
    coordinatewise, and in particular $k_h\in\Lp^2(0,T)$.
    It remains to compare the two terms in~\eqref{eq:k-h} near $T$.
    By Lemma~\ref{lem:kernel}~\ref{kernel:asymptotics}, the first term satisfies $h(0)K_H(T,r)\geq\frac12h(0)c'_H(T-r)^{H-1/2}$ for $r$ close to $T$, where $h(0)>0$ and $c'_H>0$.
    By Lemma~\ref{lem:kernel}~\ref{kernel:bound} with $\delta=T/2$, the second term is bounded in absolute value by $C\Norm{h'}_\infty(T-r)^{H+1/2}$ for $r\geq T/2$, which is of smaller order as $r\uparrow T$.
    Hence $k_h>0$ on an interval $(T-\eps,T)$ for some $\eps>0$.

    \emph{Step~3.}
    We prove~\ref{fbm:statistic} and~\ref{fbm:statistic-sde}.
    The observations $O_u\defined X_u/c(u)=X_0+G_u$, $u\in[0,T]$, are linear-Gaussian in the sense of Definition~\ref{def:observations} with $a_u=1$.
    Indeed, the coordinates of $B$ are independent fBMs, and $G$ is a linear functional of $B$.
    With $I_t\defined[t,T]$ we have $\Info_{I_t}=\Info^{\mathrm{A}}_t$, the sets $I_t$ are nested, and $I_T=\{T\}$ sees the data.
    In the construction of the proof of Proposition~\ref{prop:sufficiency} we take $i_0=T$, so that $\zeta=G^{(1)}_T$ and the signal-free combinations are $S_u=G_u-G_T$.
    By Lemma~\ref{lem:gaussian-hilbert}~\ref{gh:projection}, applied coordinatewise, the projection of $\zeta$ yields the conditional expectation in~\ref{fbm:statistic}, so that $m_{I_t}=m_t$.
    Proposition~\ref{prop:sufficiency}~\ref{suff:blue} thus shows that $m_t$ is the best linear unbiased estimator of $X_0$ from the path of $X$ on $[t,T]$.

    To compute $m_t$, recall that $\sigma(G_u-G_T:u\in[t,T])=\Filt^{\rev{W}}_{T-t}$ by Step~1 and that $G_T=-\rev{G}_T=-\int_0^Tk_h\de\rev{W}$ by Step~2.
    Hence
    \[
        \E\bigl[G_T\given G_u-G_T,\ u\in[t,T]\bigr]
        =-\E\Bigl[\int_0^Tk_h\de\rev{W}\given[\Big]\Filt^{\rev{W}}_{T-t}\Bigr]
        =-\int_0^{T-t}k_h\de\rev{W},
    \]
    which is the formula for $\rev{M}$ in~\ref{fbm:statistic-sde}, and it follows that
    \[
        m_t-X_0=G_T+\int_0^{T-t}k_h\de\rev{W}=-\int_{T-t}^Tk_h\de\rev{W}.
    \]
    By the It\^{o} isometry, the variance per coordinate of this error is $v(t)=v_{I_t}=\int_{T-t}^Tk_h(r)^2\de r=\vA(t)$.
    This function is positive by Step~2, absolutely continuous with derivative $k_h(T-t)^2$, tends to zero as $t\downarrow0$, and equals $\Var[G^{(1)}_T]$ for $t=T$.
    The bound $\vA(t)\leq\Var[G^{(1)}_t]$ holds because $X_t/c(t)=X_0+G_t$ is itself an unbiased linear statistic.
    If $h'\geq0$, both terms in~\eqref{eq:k-h} are nonnegative, and the first one is positive by Lemma~\ref{lem:kernel}~\ref{kernel:positive}, so that $k_h>0$ on $(0,T)$.
    In particular, $\vA$ satisfies the hypotheses of Theorem~\ref{thm:twin}~\ref{twin:sde}, and the right-hand side in~\ref{fbm:statistic-sde} is continuous in $s$, so it is the continuous modification used there.

    \emph{Step~4.}
    We prove~\ref{fbm:path-sde}.
    For $s\in[0,T)$, let $\rev\Info_s\defined\Info^{\mathrm{A}}_{T-s}$ and $\FiltAlt_s\defined\rev\Info_s\vee\sigma(X_0)$.
    With $\Phi\defined\int_0^Tk_h\de\rev{W}$, we have $X_T/c(T)=X_0+G_T=X_0-\Phi$, so~\eqref{eq:info-a} gives $\FiltAlt_s=\sigma(X_0,\Phi)\vee\Filt^{\rev{W}}_s$.
    The Brownian motion $\rev{W}$ is a functional of $B$ and hence independent of $X_0$.
    Lemma~\ref{lem:enlargement} in the appendix, applied with $U=X_0$, $k=k_h$, and $\rho=\rev{v}$, therefore shows that $\rev{W}-\int_0^{\argDot}b_r\de r$ is a $\FiltAlt$-Brownian motion on $[0,T)$, where
    \[
        b_r\defined\frac{k_h(r)}{\rev{v}(r)}\,R_r
        \qquad\text{with}\qquad
        R_r\defined\int_r^Tk_h\de\rev{W}=X_0-\rev{M}_r.
    \]
    Fix $S<T$.
    On $[0,S]$, the drift $b$ is $\FiltAlt$-progressively measurable, and since $\E[\Abs{R_r}]\leq\sqrt{d\,\rev{v}(r)}$ and $\rev{v}$ is nonincreasing, it satisfies
    \[
        \E\int_0^S\Abs{b_r}\de r\leq\sqrt{d}\,\rev{v}(S)^{-1/2}\int_0^S\Abs{k_h(r)}\de r<\infty.
    \]
    Moreover, $\rev{W}$ is $\rev\Info$-adapted by~\eqref{eq:info-a}, and $R_r=-e_{T-r}$ in the notation of Proposition~\ref{prop:sufficiency}.
    Part~\ref{suff:posterior} of that proposition with $J=I_{T-r}$ therefore gives
    \[
        \E\bigl[b_r\given\rev\Info_r\bigr]
        =-\frac{k_h(r)}{\rev{v}(r)}\,\E\bigl[e_{T-r}\given\Info^{\mathrm{A}}_{T-r}\bigr]
        =k_h(r)\,\Score_{\rev{v}(r)}\bigl(\rev{M}_r\bigr),
        \qquad r\in[0,S].
    \]
    Lemma~\ref{lem:innovation} in the appendix now shows that $\widehat{W}$ is a continuous $\rev\Info$-martingale on $[0,S]$.
    Since $\widehat{W}-\rev{W}$ has finite variation on $[0,S]$, we have $\langle\widehat{W}^{(i)},\widehat{W}^{(j)}\rangle_s=\delta_{ij}s$, and L\'{e}vy's characterization shows that $\widehat{W}$ is an $\rev\Info$-Brownian motion on $[0,T)$.
    As a Brownian motion on $[0,T)$, it has an almost sure limit at $T$, which gives the extension.

    It remains to prove the representation of the reversed path.
    For $s<T$ we have $\rev{X}_s=X_0+G_{T-s}=X_T/c(T)+\rev{G}_s$, so the first identity is~\eqref{eq:reversed-noise}.
    For the second identity and for the martingale part of $\rev{M}$, we show that every $f\in\Lp^2(0,s)$ satisfies
    \[
        \int_0^sf(r)\de\rev{W}_r
        =\int_0^sf(r)\de\widehat{W}_r+\int_0^sf(r)k_h(r)\,\Score_{\rev{v}(r)}\bigl(\rev{M}_r\bigr)\de r.
    \]
    The last integral converges absolutely.
    Indeed, Proposition~\ref{prop:sufficiency}~\ref{suff:posterior} gives $\E\Abs{\Score_{\rev{v}(r)}(\rev{M}_r)}\leq\sqrt{d/\rev{v}(r)}$, and the Cauchy--Schwarz inequality yields
    \[
        \E\int_0^s\Abs{f(r)k_h(r)}\,\Abs[\big]{\Score_{\rev{v}(r)}\bigl(\rev{M}_r\bigr)}\de r
        \leq\sqrt{d}\,\rev{v}(s)^{-1/2}\Norm{f}_{\Lp^2}\Norm{k_h}_{\Lp^2}.
    \]
    For step functions $f$, the identity is the definition of $\widehat{W}$.
    For general $f$, we choose step functions $f_n\to f$ in $\Lp^2(0,s)$.
    Then the Wiener integrals with respect to $\rev{W}$ and $\widehat{W}$ converge in $\Lp^2(\Prob)$ by the It\^{o} isometry, and the drift integrals converge in $\Lp^1(\Prob)$ by the bound above with $f$ replaced by $f_n-f$.
    With $f=K_H(s,\argDot)$, the identity gives the representation of $\rev{B}_s$.
    With $f=k_h$, the identity and~\ref{fbm:statistic-sde} show that
    \[
        \rev{M}_s-\rev{M}_0=\int_0^sk_h\de\widehat{W}+\int_0^sk_h(r)^2\,\Score_{\rev{v}(r)}\bigl(\rev{M}_r\bigr)\de r,
        \qquad s\in[0,T),
    \]
    so the martingale part of $\rev{M}$ is $\int_0^{\argDot}k_h\de\widehat{W}$, and the proof is complete.
\end{proof}

For $H=\frac12$ we have $k_h=h(T-\argDot)$ and $m_t=X_t/c(t)$, and Theorem~\ref{thm:fbm} is the classical time reversal of a diffusion model.
For pure fBM noising, the schedule is explicit and can be compared with the variance $t^{2H}$ of the noise in $X_t$ alone.

\begin{proposition}[Pure fBM Noising]
    \label{prop:pure-fbm}
    Let $\mu=0$ and $g=1$, so that $X_t=X_0+B_t$.
    Then the schedule of Theorem~\ref{thm:fbm} is given by
    \[
        \vA(t)=\int_{T-t}^TK_H(T,r)^2\de r,\qquad t\in[0,T],
    \]
    it is strictly increasing with $\vA(T)=T^{2H}$, and the following hold.
    \begin{enumerate}
        \item \label{pure:strict}
              If $H=\frac12$, then $\vA(t)=t$ for all $t\in[0,T]$, and if $H\neq\frac12$, then $\vA(t)<t^{2H}$ for all $t\in(0,T)$.
        \item \label{pure:asymptotics}
              As $t\downarrow0$, the ratio $\vA(t)/t^{2H}$ converges to
              \[
                  \kappa_H\defined\frac{\Gamma(\frac32-H)}{\Gamma(H+\frac12)\,\Gamma(2-2H)}.
              \]
    \end{enumerate}
    The Brownian twin is the variance-exploding model $\de\Twin_t=K_H(T,T-t)\de\widetilde{W}_t$ started at $X_0$, whose law at time $T$ is $\pi\conv\Nor(0,T^{2H}\Id)$, the same as that of $X_T$.
\end{proposition}

\begin{proof}
    Here $c\equiv1$ and $h\equiv1$, so $k_h=K_H(T,\argDot)$, and the formula for $\vA$ is Theorem~\ref{thm:fbm}~\ref{fbm:statistic}.
    Since $K_H>0$ by Lemma~\ref{lem:kernel}~\ref{kernel:positive}, the schedule is strictly increasing, and Lemma~\ref{lem:molchan-golosov} gives $\vA(T)=\int_0^TK_H(T,r)^2\de r=T^{2H}$.
    The statement on the twin follows from~\eqref{eq:twin-sde} with $\mu=0$ and $\vA'(t)=K_H(T,T-t)^2$.

    We prove~\ref{pure:strict}.
    The statistic $X_t=X_0+B_t$ is unbiased with error variance $t^{2H}$ per coordinate, so $\vA(t)\leq t^{2H}$.
    By~\eqref{eq:pythagoras} with $i_0=T$, equality holds if and only if $B^{(1)}_t$ is orthogonal to all signal-free combinations, that is, to $B^{(1)}_u-B^{(1)}_T$ for all $u\in[t,T]$.
    These covariances are
    \[
        \Cov\bigl[B^{(1)}_t,B^{(1)}_u-B^{(1)}_T\bigr]=R_H(t,u)-R_H(t,T)=\tfrac12\bigl(f(u)-f(T)\bigr),
        \quad f(u)\defined u^{2H}-(u-t)^{2H}.
    \]
    For $u>t$, the derivative $f'(u)=2H(u^{2H-1}-(u-t)^{2H-1})$ vanishes if and only if $H=\frac12$.
    So for $H\neq\frac12$ and $t<T$, the function $f$ is strictly monotone on $[t,T]$, the covariance is nonzero for $u\in[t,T)$, and $\vA(t)<t^{2H}$.
    For $H=\frac12$, the covariance vanishes and $\vA(t)=t$.

    We prove~\ref{pure:asymptotics}.
    By Lemma~\ref{lem:kernel}~\ref{kernel:asymptotics}, we have $K_H(T,r)^2/(T-r)^{2H-1}\to(c'_H)^2$ as $r\uparrow T$, and therefore
    \[
        \vA(t)=\int_0^tK_H(T,T-q)^2\de q\sim(c'_H)^2\int_0^tq^{2H-1}\de q=\frac{(c'_H)^2}{2H}\,t^{2H}
        \qquad\text{as }t\downarrow0.
    \]
    Since $(c'_H)^2/(2H)=\kappa_H$ by Lemma~\ref{lem:kernel}~\ref{kernel:kappa}, the proof is complete.
\end{proof}

\section{Generative Fractional Diffusion Models}
\label{sec:gfdm}

We now turn to the generative fractional diffusion models (GFDM) of \citet{nobis2025generative.local}, in which the noise is a weighted sum of Ornstein--Uhlenbeck processes driven by one Brownian motion.
Here the reverse-time process runs on an augmented state, which is Markovian, so the question is how much of the augmented state is needed.
We show that a single linear statistic of it is sufficient, identify the schedule of the Brownian twin together with its derivative, and characterize when this schedule degenerates.

Let $K\in\N$, let $0<\gamma_1<\dots<\gamma_K$ be speeds and $\omega=(\omega_1,\dots,\omega_K)\in\R^K$ weights, and let $B$ be a standard Brownian motion in $\R^d$ independent of $X_0$.
Following \citet{nobis2025generative.local}, we define the Ornstein--Uhlenbeck processes $Y^1,\dots,Y^K$ and the Markovian approximation $\widehat{B}$ of fractional Brownian motion of type~II by
\[
    Y^k_t\defined\int_0^te^{-\gamma_k(t-r)}\de B_r
    \qquad\text{and}\qquad
    \widehat{B}_t\defined\sum_{k=1}^K\omega_kY^k_t,
    \qquad t\in[0,T].
\]
In this setting, the noise in~\eqref{eq:forward} is $N=\widehat{B}$, with $\mu$, $g$, $c$, and $h=g/c$ as in Section~\ref{sec:setting}, so that
\[
    X_t=c(t)\bigl(X_0+\Gamma_t\bigr)
    \qquad\text{with}\qquad
    \Gamma_t\defined\int_0^th(r)\de\widehat{B}_r,
    \qquad t\in[0,T],
\]
where the integral is an It\^{o} integral.
The \emph{augmented state} $Z_t\defined(X_t,Y^1_t,\dots,Y^K_t)\in(\R^d)^{1+K}$ solves a linear stochastic differential equation with continuous coefficients driven by $B$, see equations~(33)--(36) in \cite{nobis2025generative.local}.
The reverse-time process of GFDM runs on the augmented state, so its information at time $t$ is
\[
    \Info^{\mathrm{B}}_t\defined\sigma\bigl(Z_u : u\in[t,T]\bigr),
    \qquad t\in[0,T].
\]

Given $X_0$, the augmented state $Z_t$ is Gaussian with mean $(c(t)X_0,0,\dots,0)$ and covariance $\Sigma_t\otimes\Id$, where the matrix $\Sigma_t\in\R^{(1+K)\times(1+K)}$ does not depend on $X_0$.
We denote its blocks by $\Sigma^{xx}_t\in\R$, $\Sigma^{xy}_t=(\Sigma^{yx}_t)^\tp\in\R^{1\times K}$, and $\Sigma^{yy}_t\in\R^{K\times K}$.
For $t>0$ the block $\Sigma^{yy}_t$ is positive definite, since $a^\tp\Sigma^{yy}_ta=\int_0^t(\sum_ka_ke^{-\gamma_k(t-r)})^2\de r$ for $a\in\R^K$ and exponentials with distinct rates are linearly independent on every interval.
We may therefore define
\begin{equation}
    \label{eq:gfdm-eta-v}
    \eta_t\defined\bigl(\Sigma^{yy}_t\bigr)^{-1}\Sigma^{yx}_t\in\R^K
    \qquad\text{and}\qquad
    \vB(t)\defined\frac{\Sigma^{xx}_t-\Sigma^{xy}_t(\Sigma^{yy}_t)^{-1}\Sigma^{yx}_t}{c(t)^2},
    \qquad t\in(0,T],
\end{equation}
together with $\vB(0)\defined0$, and we write $\eta_t\cdot Y_t\defined\sum_k\eta^k_tY^k_t$.
In the notation of \citet{nobis2025generative.local}, the $\eta^k_t$ are the coefficients of their equation~(67), $\Sigma^{xx}_t=c(t)^2\sigma^2(t)$, and $\Sigma^{xy}_t(\Sigma^{yy}_t)^{-1}\Sigma^{yx}_t=\tau^2_t$, so that $c(t)^2\vB(t)=c(t)^2\sigma^2(t)-\tau^2_t$.

Finally, we describe the augmented path as a family of observations in the sense of Definition~\ref{def:observations}.
It consists of
\[
    O_{(u,0)}\defined\frac{X_u}{c(u)}=X_0+\Gamma_u
    \qquad\text{and}\qquad
    O_{(u,k)}\defined Y^k_u,\quad k=1,\dots,K,
\]
indexed by $(u,k)\in I\defined[0,T]\times\{0,\dots,K\}$.
When the reverse-time process has reached time $t$, it has observed those indexed by $I_t\defined[t,T]\times\{0,\dots,K\}$.

\begin{assumption}
    \label{ass:gfdm-nondegenerate}
    We assume that $\vB(t)>0$ for every $t\in(0,T]$.
\end{assumption}

The assumption is needed for the statistic of Theorem~\ref{thm:gfdm} below to have a density, so that its posterior and its score are defined.
It is not automatic, and Proposition~\ref{prop:gfdm-degenerate} below shows that it can fail.

Since the augmented state is Markovian, one expects the augmented path on $[t,T]$ to carry no more information about the data than the current state $Z_t$.
The question is then how much of $Z_t$ is needed.
The following theorem shows that the answer is a single vector, the residual of the regression of $X_t$ on $Y_t$, and that this vector is the data plus independent Gaussian noise of variance $\vB(t)$.

\begin{theorem}[Reduction of GFDM]
    \label{thm:gfdm}
    Suppose that Assumption~\ref{ass:gfdm-nondegenerate} holds, let $t\in(0,T]$, and let
    \[
        \tilde{x}_t\defined X_t-\eta_t\cdot Y_t
        \qquad\text{and}\qquad
        m_t\defined\frac{\tilde{x}_t}{c(t)}.
    \]
    Then the following hold.
    \begin{enumerate}
        \item \label{gfdm:fixed-time}
              The statistic satisfies $m_t=X_0+e_t$ with $e_t\sim\Nor(0,\vB(t)\Id)$, and $X_0$, $e_t$, and $Y_t$ are independent.
              The density of $Z_t$ factorizes as
              \[
                  p_t(x,y)=q_t(y)\,c(t)^{-d}\,\Smoothed_{\vB(t)}\Bigl(\frac{x-\eta_t\cdot y}{c(t)}\Bigr),
                  \qquad (x,y)\in\R^d\times(\R^d)^K,
              \]
              where $q_t$ is the Gaussian density of $Y_t$, and consequently, for $k=1,\dots,K$,
              \begin{align}
                  \nabla_x\log p_t(x,y)     & =\frac{1}{c(t)}\,\Score_{\vB(t)}\Bigl(\frac{x-\eta_t\cdot y}{c(t)}\Bigr), \\
                  \nabla_{y^k}\log p_t(x,y) & =\nabla_{y^k}\log q_t(y)-\eta^k_t\,\nabla_x\log p_t(x,y).
              \end{align}
        \item \label{gfdm:markov}
              The conditional law of $X_0$ given the augmented path on $[t,T]$ depends only on $m_t$, namely
              \[
                  \Prob\bigl(X_0\in\de x\given\Info^{\mathrm{B}}_t\bigr)
                  =\Prob\bigl(X_0\in\de x\given Z_t\bigr)
                  =\pi_{\vB(t),m_t}(\de x).
              \]
              Moreover, $m_t$ is the best linear unbiased estimator of $X_0$ from the augmented path on $[t,T]$, that is, $m_t=m_{I_t}$ and $\vB(t)=v_{I_t}$ in the notation of Theorem~\ref{thm:twin}.
        \item \label{gfdm:schedule}
              The function $\vB$ is nondecreasing and absolutely continuous on $[0,T]$, continuously differentiable on $(0,T]$, and satisfies $\vB(t)\leq\Sigma^{xx}_t/c(t)^2$.
    \end{enumerate}
    In particular, Theorem~\ref{thm:twin} and Corollary~\ref{cor:score-matching} apply with $v=\vB$ and $\Info_t=\Info^{\mathrm{B}}_t$.
\end{theorem}

\begin{proof}
    The last claim follows from parts~\ref{gfdm:markov} and~\ref{gfdm:schedule} and Assumption~\ref{ass:gfdm-nondegenerate}, together with Step~1 below, which puts the augmented path into the situation of Theorem~\ref{thm:twin}; it remains to prove~\ref{gfdm:fixed-time}--\ref{gfdm:schedule}.
    We first split the augmented path on $[t,T]$ into the current state $Z_t$ and innovations which are independent of the past, and then treat the three parts in turn.

    \emph{Step~1.} We show that the observations are linear-Gaussian and split them at time $t$.
    Their noises are $G_{(u,0)}=\Gamma_u$ and $G_{(u,k)}=Y^k_u$, which are linear functionals of $B$.
    Since $B$ is independent of $X_0$ and has independent coordinates, Definition~\ref{def:observations} holds with $a_{(u,0)}=1$ and $a_{(u,k)}=0$ for $k\geq1$.
    Moreover, the sets $I_t$ decrease in $t$, we have $\Info_{I_t}=\Info^{\mathrm{B}}_t$, and $I_T$ sees the data.
    Now fix $t$ and let $u\in[t,T]$.
    From $\de Y^k_r=-\gamma_kY^k_r\de r+\de B_r$ and $\de\Gamma_r=h(r)\sum_k\omega_k\de Y^k_r$ we obtain
    \begin{equation}
        \label{eq:gfdm-markov}
        Y^k_u=e^{-\gamma_k(u-t)}Y^k_t+\varpi^k_{u,t}
        \qquad\text{and}\qquad
        \frac{X_u}{c(u)}=\frac{X_t}{c(t)}+\alpha_{u,t}\cdot Y_t+\varpi^0_{u,t},
    \end{equation}
    where $\varpi_{u,t}=(\varpi^0_{u,t},\dots,\varpi^K_{u,t})$ is a linear functional of $(B_r-B_t)_{r\in[t,u]}$.
    The vector $\alpha_{u,t}\in\R^K$ is deterministic with entries $\alpha^k_{u,t}=-\omega_k\gamma_k\int_t^uh(r)e^{-\gamma_k(r-t)}\de r$.
    Let $(\FiltAlt_t)_{t\in[0,T]}$ be the augmented natural filtration of $B$.
    Then $(\varpi_{u,t})_{u\in[t,T]}$ is independent of $\FiltAlt_t\vee\sigma(X_0)$.
    By~\eqref{eq:gfdm-markov}, each $Z_u$ is a function of $Z_t$ and $\varpi_{u,t}$, and conversely $\varpi_{u,t}$ is a function of $Z_t$ and $Z_u$, so that
    \begin{equation}
        \label{eq:gfdm-info}
        \Info^{\mathrm{B}}_t=\sigma(Z_t)\vee\sigma\bigl(\varpi_{u,t} : u\in[t,T]\bigr).
    \end{equation}

    \emph{Step~2.} We prove~\ref{gfdm:fixed-time}.
    Given $X_0$, the noise of $Z_t$ is $(c(t)\Gamma_t,Y_t)$, whose covariance is $\Sigma_t\otimes\Id$.
    Let $\eps_t\defined c(t)\Gamma_t-\eta_t\cdot Y_t$, and let $Y^{(1)}_t\in\R^K$ collect the first coordinates of $Y^1_t,\dots,Y^K_t$.
    By~\eqref{eq:gfdm-eta-v}, the first coordinate of $\eps_t$ satisfies
    \begin{align}
        \Cov\bigl[\eps^{(1)}_t,Y^{(1)}_t\bigr] & =\Sigma^{xy}_t-\eta_t^\tp\Sigma^{yy}_t=0,                                             \\
        \Var\bigl[\eps^{(1)}_t\bigr]           & =\Sigma^{xx}_t-\Sigma^{xy}_t\bigl(\Sigma^{yy}_t\bigr)^{-1}\Sigma^{yx}_t=c(t)^2\vB(t),
    \end{align}
    and the same holds in every other coordinate.
    As different coordinates are independent, $\eps_t$ and $Y_t$ are uncorrelated jointly Gaussian vectors, hence independent, and both are independent of $X_0$ as functionals of $B$.
    Since $\tilde{x}_t=c(t)X_0+\eps_t$, we obtain $m_t=X_0+e_t$ with $e_t\defined\eps_t/c(t)\sim\Nor(0,\vB(t)\Id)$, and $X_0$, $e_t$, and $Y_t$ are independent.
    By Assumption~\ref{ass:gfdm-nondegenerate}, $m_t$ therefore has the density $\Smoothed_{\vB(t)}$, and since $\Sigma^{yy}_t$ is positive definite, $Y_t$ has a Gaussian density $q_t$.
    The linear map $(m,y)\mapsto(c(t)m+\eta_t\cdot y,y)$ has Jacobian determinant $c(t)^d$ and maps $(m_t,Y_t)$ to $Z_t$, which gives the factorization of $p_t$.
    The two formulas for the scores follow from the factorization by the chain rule.

    \emph{Step~3.} We prove~\ref{gfdm:markov}.
    We first identify the conditional laws.
    By~\eqref{eq:gfdm-info} and $\sigma(Z_t)=\sigma(m_t)\vee\sigma(Y_t)$, we have
    \[
        \Info^{\mathrm{B}}_t=\sigma(m_t)\vee\mathcal{S}_t
        \qquad\text{with}\qquad
        \mathcal{S}_t\defined\sigma\bigl(Y_t,\varpi_{u,t} : u\in[t,T]\bigr).
    \]
    The variables $X_0$, $e_t$, $Y_t$, and $(\varpi_{u,t})_{u\in[t,T]}$ are independent, because $e_t$ and $Y_t$ are $\FiltAlt_t$-measurable and, by Step~2, independent of each other and of $X_0$, while $(\varpi_{u,t})_{u\in[t,T]}$ is independent of $\FiltAlt_t\vee\sigma(X_0)$ by Step~1.
    Hence $(X_0,m_t)$ is independent of $\mathcal{S}_t$, and the argument in the proof of Proposition~\ref{prop:sufficiency}~\ref{suff:posterior} shows that the conditional law of $X_0$ given $\Info^{\mathrm{B}}_t$, and likewise given $Z_t$, is the conditional law given $m_t$, which is $\pi_{\vB(t),m_t}$.

    It remains to show that $m_t$ is the best linear unbiased estimator.
    Let $\mathcal{L}_t$ be the closed linear span of the first coordinates of $Y_t$ and of $\varpi_{u,t}$, $u\in[t,T]$, which is independent of and hence orthogonal to the first coordinate $e^{(1)}_t$ of $e_t$.
    By~\eqref{eq:gfdm-markov}, the terms $X_t/c(t)$ cancel in every signal-free combination of observations in $I_t$, that is, in every $c\cdot O$ with $\sum_ic_ia_i=0$, so the first coordinate of its noise lies in $\mathcal{L}_t$.
    Now let $c\cdot O$ be unbiased.
    Then $c\cdot O-m_t$ is signal-free, so $c\cdot G^{(1)}=e^{(1)}_t+\ell$ with $\ell\in\mathcal{L}_t$, and
    \[
        \E\bigl[(c\cdot G^{(1)})^2\bigr]=\vB(t)+\E\bigl[\ell^2\bigr]\geq\vB(t).
    \]
    Since $m_t$ is itself an unbiased linear statistic, this shows that $\vB(t)=v_{I_t}$, and Proposition~\ref{prop:sufficiency}~\ref{suff:blue}, applied to the coefficients of $m_t$, gives $m_t=m_{I_t}$ almost surely.

    \emph{Step~4.} We prove~\ref{gfdm:schedule}.
    We first show that $\vB$ is nondecreasing.
    By part~\ref{gfdm:markov}, we have $\vB(t)=v_{I_t}$ for $t\in(0,T]$.
    For $0<t'\leq t$, the set $I_{t'}$ contains $I_t$, so the infimum defining $v_{I_{t'}}$ runs over more unbiased linear statistics than the one defining $v_{I_t}$, and hence $\vB(t')\leq\vB(t)$.
    Since $\vB(0)=0$ and $\vB\geq0$, the function $\vB$ is nondecreasing on $[0,T]$.

    Next, we show that $\Sigma$ is continuously differentiable on $[0,T]$.
    Let $\chi_t\in\R^{1+K}$ collect the first coordinates of the noise $(c(t)\Gamma_t,Y^1_t,\dots,Y^K_t)$ of $Z_t$, so that $\Sigma_t=\E[\chi_t\chi_t^\tp]$.
    From $\de Y^k_t=-\gamma_kY^k_t\de t+\de B_t$, $\de\Gamma_t=h(t)\sum_k\omega_k\de Y^k_t$, $c'=\mu c$, and $ch=g$, we obtain the linear equation
    \[
        \de\chi_t=F(t)\chi_t\de t+G(t)\de B^{(1)}_t,\qquad\chi_0=0,
    \]
    with the continuous coefficients
    \[
        F(t)=\begin{pmatrix}
            \mu(t) & -g(t)\omega_1\gamma_1 & \cdots & -g(t)\omega_K\gamma_K \\
            0      & -\gamma_1             &        &                       \\
            \vdots &                       & \ddots &                       \\
            0      &                       &        & -\gamma_K
        \end{pmatrix}
        \qquad\text{and}\qquad
        G(t)=\begin{pmatrix}g(t)\sum_k\omega_k\\1\\\vdots\\1\end{pmatrix}.
    \]
    It\^{o}'s product rule gives $\de(\chi_t\chi_t^\tp)=(F\chi_t\chi_t^\tp+\chi_t\chi_t^\tp F^\tp+GG^\tp)(t)\de t+(G\chi_t^\tp+\chi_tG^\tp)(t)\de B^{(1)}_t$.
    The stochastic integral has mean zero because its integrand is square integrable, so taking expectations yields the Lyapunov equation
    \[
        \Sigma_t=\int_0^t\bigl(F(r)\Sigma_r+\Sigma_rF(r)^\tp+G(r)G(r)^\tp\bigr)\de r,\qquad t\in[0,T].
    \]
    The integrand is bounded, so $\Sigma$ is continuous with $\Sigma_0=0$.
    Hence the integrand is continuous as well, and $\Sigma$ is continuously differentiable on $[0,T]$.

    We now turn to the regularity of $\vB$.
    For $t>0$, the block $\Sigma^{yy}_t$ is positive definite, and matrix inversion is smooth on the set of invertible matrices, so $t\mapsto(\Sigma^{yy}_t)^{-1}$ is continuously differentiable on $(0,T]$.
    Moreover, $c$ is positive and continuously differentiable, since $c'=\mu c$ with $\mu$ continuous.
    By~\eqref{eq:gfdm-eta-v}, the function $\vB$ is therefore continuously differentiable on $(0,T]$.
    Since $(\Sigma^{yy}_t)^{-1}$ is positive definite, the subtracted term in~\eqref{eq:gfdm-eta-v} is nonnegative, and therefore
    \[
        0\leq\vB(t)\leq\frac{\Sigma^{xx}_t}{c(t)^2},\qquad t\in(0,T].
    \]
    The right-hand side is the noise variance of the unbiased statistic $X_t/c(t)$ alone.
    As $\Sigma$ is continuous with $\Sigma_0=0$ and $c(0)=1$, it tends to zero as $t\downarrow0$, so $\vB$ is continuous also at $t=0$.

    Finally, we show that $\vB$ is absolutely continuous on $[0,T]$.
    As $\vB$ is nondecreasing and continuously differentiable on $(0,T]$, its derivative is nonnegative there, and the fundamental theorem of calculus gives
    \[
        \vB(t)-\vB(\eps)=\int_\eps^t\vB'(r)\de r,\qquad 0<\eps<t\leq T.
    \]
    As $\eps\downarrow0$, the left-hand side tends to $\vB(t)$ by the continuity at zero, and the right-hand side tends to $\int_0^t\vB'(r)\de r$ by monotone convergence.
    Hence $\vB'$ is integrable on $(0,T]$ and $\vB(t)=\int_0^t\vB'(r)\de r$ for every $t\in[0,T]$, that is, $\vB$ is absolutely continuous on $[0,T]$ with $\vB(0)=0$.
    This is the regularity required in Theorem~\ref{thm:twin}~\ref{twin:sde}, and the proof is complete.
\end{proof}

The vector $\tilde{x}_t=c(t)m_t$, which is the network input of Corollary~\ref{cor:score-matching}, is the network input used by GFDM.
Moreover, by part~\ref{gfdm:fixed-time}, the law of $X_t$ given $(X_0,Y_t)$ is $\Nor(c(t)X_0+\eta_t\cdot Y_t,c(t)^2\vB(t)\Id)$, so the target of the augmented score matching loss of \cite{nobis2025generative.local} is
\[
    \nabla_x\log p_{t|0}\bigl(X_t\given Y_t,X_0\bigr)=-\frac{\tilde{x}_t-c(t)X_0}{c(t)^2\vB(t)}=\Lambda_t,
\]
and this loss is the loss $\LossFrac$ of Corollary~\ref{cor:score-matching} with constant weight and uniform time.

Whether Assumption~\ref{ass:gfdm-nondegenerate} holds is decided by a conserved quantity of the augmented dynamics.

\begin{proposition}[Degeneracy]
    \label{prop:gfdm-degenerate}
    It holds that
    \begin{equation}
        \label{eq:gfdm-conserved}
        \frac{X_t}{c(t)}-h(t)\sum_{k=1}^K\omega_kY^k_t=X_0-\int_0^th'(r)\widehat{B}_r\de r,
        \qquad t\in[0,T].
    \end{equation}
    \begin{enumerate}
        \item \label{degenerate:constant}
              If $h$ is constant on $[0,t_0]$ for some $t_0>0$, then $\vB=0$ on $[0,t_0]$, and for $t\in[0,t_0]$ the data $X_0$ equal the left-hand side of~\eqref{eq:gfdm-conserved}, which is a linear function of the augmented state $Z_t$.
        \item \label{degenerate:exponential}
              Let $\mu=0$ and $g(t)=a\,e^{\lambda t}$ with $a,\lambda>0$ and $\lambda\notin\{\gamma_1,\dots,\gamma_K\}$.
              Then, for every $t\in(0,T]$,
              \[
                  \vB(t)>0
                  \qquad\text{if and only if}\qquad
                  \sum_{k=1}^K\frac{\omega_k}{\lambda-\gamma_k}\neq0.
              \]
              In particular, Assumption~\ref{ass:gfdm-nondegenerate} holds if and only if this sum is nonzero.
    \end{enumerate}
\end{proposition}

\begin{proof}
    We first prove~\eqref{eq:gfdm-conserved}.
    Since $h$ is deterministic and continuously differentiable, It\^{o}'s product rule gives
    \[
        \Gamma_t=\int_0^th(r)\de\widehat{B}_r=h(t)\widehat{B}_t-\int_0^th'(r)\widehat{B}_r\de r,
    \]
    and~\eqref{eq:gfdm-conserved} follows from $X_t/c(t)=X_0+\Gamma_t$ and $\widehat{B}_t=\sum_k\omega_kY^k_t$.

    We prove~\ref{degenerate:constant}.
    If $h'=0$ on $[0,t_0]$, then the right-hand side of~\eqref{eq:gfdm-conserved} equals $X_0$ for $t\leq t_0$.
    Hence $X_0$ is an unbiased linear statistic of $Z_t$ without error, and $\vB(t)=0$.

    We prove~\ref{degenerate:exponential}.
    Here $c\equiv1$ and $h=g$.
    By Step~2 of the proof of Theorem~\ref{thm:gfdm}, $\vB(t)$ is the variance of the residual of $\Gamma^{(1)}_t$ after regression on $Y^{(1)}_t$, so $\vB(t)=0$ if and only if $\Gamma^{(1)}_t$ lies in the span of the first coordinates $Y^{(1),1}_t,\dots,Y^{(1),K}_t$ of $Y^1_t,\dots,Y^K_t$.
    By~\eqref{eq:gfdm-conserved}, this is the case if and only if $J_t\defined\int_0^th'(r)\widehat{B}^{(1)}_r\de r$ lies in this span.
    We write $Y^{(1),k}_r=\int_0^re^{-\gamma_k(r-u)}\de B^{(1)}_u$ and exchange the integrals, which is justified because the Wiener integral is an isometry from $\Lp^2(0,t)$ into $\Lp^2(\Prob)$ and hence commutes with the Bochner integral of the continuous $\Lp^2(\Prob)$-valued map $r\mapsto h'(r)\widehat{B}^{(1)}_r$.
    With $h'(r)=a\lambda e^{\lambda r}$, a direct calculation gives $J_t=\int_0^t\psi(u)\de B^{(1)}_u$ with
    \[
        \psi(u)=\int_u^ta\lambda e^{\lambda r}\sum_{k=1}^K\omega_ke^{-\gamma_k(r-u)}\de r
        =a\lambda\sum_{k=1}^K\frac{\omega_k}{\lambda-\gamma_k}\bigl(e^{\lambda t}e^{-\gamma_k(t-u)}-e^{\lambda u}\bigr),
        \quad u\in[0,t].
    \]
    The elements of the span are the integrals $\int_0^t\sum_kb_ke^{-\gamma_k(t-u)}\de B^{(1)}_u$ with $b\in\R^K$, so by the It\^{o} isometry, $J_t$ lies in the span if and only if $\psi$ agrees almost everywhere on $[0,t]$ with a function $u\mapsto\sum_kb_ke^{\gamma_ku}$.
    Since $\lambda,\gamma_1,\dots,\gamma_K$ are distinct, the functions $e^{\lambda u},e^{\gamma_1u},\dots,e^{\gamma_Ku}$ are linearly independent on $[0,t]$.
    Hence $J_t$ lies in the span if and only if the coefficient $-a\lambda\sum_k\omega_k/(\lambda-\gamma_k)$ of $e^{\lambda u}$ in $\psi$ vanishes, that is, as $a\lambda>0$, if and only if $\sum_k\omega_k/(\lambda-\gamma_k)=0$.
\end{proof}

For the variance-exploding choice of \citet{nobis2025generative.local}, $c\equiv1$ and $g$ is exponential, so part~\ref{degenerate:exponential} decides Assumption~\ref{ass:gfdm-nondegenerate} whenever $\lambda\notin\{\gamma_1,\dots,\gamma_K\}$.
For the variance-preserving choice, $h=\sqrt{\beta}/c$ is strictly increasing, so part~\ref{degenerate:constant} does not force degeneracy.

Under Assumption~\ref{ass:gfdm-nondegenerate}, the derivative of the schedule is explicit, and with it the diffusion coefficient of the Brownian twin.

\begin{proposition}[Diffusion Coefficient of the Twin]
    \label{prop:gfdm-derivative}
    Suppose that Assumption~\ref{ass:gfdm-nondegenerate} holds.
    Then, for $t\in(0,T]$,
    \[
        \vB'(t)=\theta_t^2
        \qquad\text{with}\qquad
        \theta_t\defined\frac{1}{c(t)}\Bigl(g(t)\sum_{k=1}^K\omega_k-\sum_{k=1}^K\eta^k_t\Bigr),
    \]
    and the Brownian model $\Twin$ with schedule $(c,\vB)$ solves
    \[
        \de\Twin_t=\mu(t)\Twin_t\de t+c(t)\Abs{\theta_t}\de\widetilde{W}_t
    \]
    for a Brownian motion $\widetilde{W}$.
\end{proposition}

\begin{proof}
    We compute the quadratic variation of the first coordinate of $e_t=m_t-X_0$ in two ways.
    Fix $\eps\in(0,T)$.
    By Step~2 of the proof of Theorem~\ref{thm:gfdm}, we have $e_t=\Gamma_t-(\eta_t\cdot Y_t)/c(t)$ for $t\in(0,T]$, and $t\mapsto\eta_t/c(t)$ is continuously differentiable on $(0,T]$.
    Using $\de\Gamma_t=h(t)\sum_k\omega_k(-\gamma_kY^k_t\de t+\de B_t)$, $\de Y^k_t=-\gamma_kY^k_t\de t+\de B_t$, and It\^{o}'s product rule, we see that $e$ is a continuous semimartingale on $[\eps,T]$ whose martingale part is
    \[
        \int_\eps^t\Bigl(h(r)\sum_{k=1}^K\omega_k-\frac{1}{c(r)}\sum_{k=1}^K\eta^k_r\Bigr)\de B_r
        =\int_\eps^t\theta_r\de B_r,
        \qquad t\in[\eps,T],
    \]
    where we have used $h=g/c$.
    Hence
    \[
        \bigl[e^{(1)}\bigr]_t-\bigl[e^{(1)}\bigr]_\eps=\int_\eps^t\theta_r^2\de r,
        \qquad t\in[\eps,T].
    \]
    On the other hand, $m_t=m_{I_t}$ by Theorem~\ref{thm:gfdm}~\ref{gfdm:markov}, so by Theorem~\ref{thm:twin}~\ref{twin:law}, $e^{(1)}$ is a continuous centered Gaussian process with covariance $\vB(t\wedge t')$.
    The computation of the quadratic variation in the proof of Lemma~\ref{lem:bridge} uses only this covariance, and applied to the time reversal $s\mapsto e^{(1)}_{T-s}$ with $\rho(s)=\vB(T-s)$, it gives
    \[
        \bigl[e^{(1)}\bigr]_t-\bigl[e^{(1)}\bigr]_\eps=\vB(t)-\vB(\eps),
        \qquad t\in[\eps,T].
    \]
    As $\theta$ is continuous on $(0,T]$ and $\eps$ is arbitrary, the two expressions yield $\vB'=\theta^2$ on $(0,T]$.
    Finally, the equation for $\Twin$ is~\eqref{eq:twin-sde} with $c\sqrt{\vB'}=c\Abs{\theta}$.
\end{proof}

It remains to discuss the law from which the reverse-time dynamics is started.
In practice, the law of the state at time $T$ is replaced by a centered Gaussian law, namely by the law $\Nor(0,\Sigma_T\otimes\Id)$ of $Z_T$ given $X_0=0$ for the augmented state and by $\Nor(0,c(T)^2\vB(T)\Id)$ for the Brownian twin.
We describe such a start by a change of measure whose density is a function of the state at time $T$.
It changes the law of this state and leaves the conditional law of the path given this state unchanged, so the new measure is the law of the exact reverse-time dynamics started from the new law.
The following result shows that the two Gaussian starts are the same approximation, so that GFDM and its Brownian twin generate the same distribution also in this case.

\begin{proposition}[Gaussian Initial Law]
    \label{prop:gfdm-initial}
    Suppose that Assumption~\ref{ass:gfdm-nondegenerate} holds, let $\Twin$ be the Brownian diffusion model with schedule $(c,\vB)$, and let $\psi\defined\GaussDensity_{\vB(T)}/\Smoothed_{\vB(T)}$.
    Then the following hold.
    \begin{enumerate}
        \item \label{initial:density}
              The law $\Nor(0,\Sigma_T\otimes\Id)$ is absolutely continuous with respect to the law of $Z_T$, and its density depends on the augmented state only through the statistic, namely
              \[
                  \frac{\de\Nor(0,\Sigma_T\otimes\Id)}{\de\Law(Z_T)}(x,y)=\psi\Bigl(\frac{x-\eta_T\cdot y}{c(T)}\Bigr),
                  \qquad (x,y)\in\R^d\times(\R^d)^K.
              \]
              In particular, $\widetilde{\Prob}\defined\psi(m_T)\,\Prob$ is a probability measure under which $Z_T\sim\Nor(0,\Sigma_T\otimes\Id)$.
        \item \label{initial:law}
              Similarly, the measure $\widetilde{\Prob}^{\mathrm{B}}\defined\psi(\Twin_T/c(T))\,\Prob$ is a probability measure, and the state $\Twin_T$ of the twin has the law $\Nor(0,c(T)^2\vB(T)\Id)$ under it.
              Moreover, the processes
              \[
                  \bigl(X_0,(m_t)_{t\in(0,T]}\bigr)\text{ under }\widetilde{\Prob}
                  \qquad\text{and}\qquad
                  \bigl(X_0,(\Twin_t/c(t))_{t\in(0,T]}\bigr)\text{ under }\widetilde{\Prob}^{\mathrm{B}}
              \]
              have the same finite-dimensional distributions.
              In particular, the generated sample $X_0$ has the same law under both measures, which is given by
              \[
                  \widetilde{\Prob}(X_0\in A)=\int_{\R^d}\pi_{\vB(T),m}(A)\,\GaussDensity_{\vB(T)}(m)\de m
              \]
              for every Borel set $A\subseteq\R^d$.
    \end{enumerate}
\end{proposition}

\begin{proof}
    We write $v_T\defined\vB(T)$ and $m(x,y)\defined(x-\eta_T\cdot y)/c(T)$, so that $m_T=m(Z_T)$.

    We prove~\ref{initial:density}.
    By Theorem~\ref{thm:gfdm}~\ref{gfdm:fixed-time}, the state $Z_T$ has the strictly positive density $p_T(x,y)=q_T(y)\,c(T)^{-d}\,\Smoothed_{v_T}(m(x,y))$.
    The law $\Nor(0,\Sigma_T\otimes\Id)$ is the law of the noise $(c(T)\Gamma_T,Y_T)$ of $Z_T$.
    By Step~2 of the proof of Theorem~\ref{thm:gfdm}, we have $c(T)\Gamma_T=\eps_T+\eta_T\cdot Y_T$, where $\eps_T/c(T)\sim\Nor(0,v_T\Id)$ is independent of $Y_T$.
    The change of variables used there therefore shows that $\Nor(0,\Sigma_T\otimes\Id)$ has the density $q_T(y)\,c(T)^{-d}\,\GaussDensity_{v_T}(m(x,y))$.
    The quotient of the two densities is $\psi(m(x,y))$, which proves the first claim, and the second claim follows from $m_T=m(Z_T)$.

    We prove~\ref{initial:law}.
    Under $\Prob$, the random vector $\Twin_T/c(T)=X_0+W_{v_T}$ has the density $\Smoothed_{v_T}$.
    Hence $\widetilde{\Prob}^{\mathrm{B}}$ is a probability measure under which this vector has the density $\psi\,\Smoothed_{v_T}=\GaussDensity_{v_T}$.
    By Theorem~\ref{thm:gfdm}~\ref{gfdm:markov} and Theorem~\ref{thm:twin}~\ref{twin:law}, the processes $(X_0,(m_t)_{t\in(0,T]})$ and $(X_0,(\Twin_t/c(t))_{t\in(0,T]})$ have the same finite-dimensional distributions under $\Prob$.
    Since the densities of $\widetilde{\Prob}$ and $\widetilde{\Prob}^{\mathrm{B}}$ are the same function $\psi$ of the values of these processes at time $T$, their finite-dimensional distributions coincide under the new measures as well.
    Finally, the conditional law of $X_0$ given $m_T$ under $\Prob$ is $\pi_{v_T,m_T}$, and $m_T$ has the density $\Smoothed_{v_T}$ by Theorem~\ref{thm:gfdm}~\ref{gfdm:fixed-time}.
    For every Borel set $A\subseteq\R^d$, we therefore obtain
    \begin{align*}
        \widetilde{\Prob}(X_0\in A)
        =\E\bigl[\psi(m_T)\,\pi_{v_T,m_T}(A)\bigr]
         & =\int_{\R^d}\pi_{v_T,m}(A)\,\psi(m)\,\Smoothed_{v_T}(m)\de m \\
         & =\int_{\R^d}\pi_{v_T,m}(A)\,\GaussDensity_{v_T}(m)\de m,
    \end{align*}
    and the proof is complete.
\end{proof}

The proposition shows which noise level matters for the Gaussian start.

\begin{remark}
    \label{rem:gfdm-terminal}
    The reverse-time dynamics of GFDM is started from the centered Gaussian law with covariance $\Sigma_T\otimes\Id$, that is, from the law of $Z_T$ given $X_0=0$.
    By Proposition~\ref{prop:gfdm-initial}, this replaces the true law $\pi\conv\Nor(0,\vB(T)\Id)$ of the statistic $m_T$ by $\Nor(0,\vB(T)\Id)$, exactly as the Gaussian start of the Brownian twin does.
    The quality of this approximation is therefore governed by $\vB(T)$ relative to the spread of the data, and not by the marginal noise variance $\Sigma^{xx}_T/c(T)^2$ of $X_T/c(T)$.
    By Theorem~\ref{thm:gfdm}~\ref{gfdm:schedule}, $\vB(T)$ is the smaller of these two noise levels, and by~\eqref{eq:gfdm-eta-v} strictly so unless $\Sigma^{xy}_T=0$.
    Separately, by~\eqref{eq:gfdm-conserved} the augmented state determines $X_0$ up to the error $\int_0^Th'(r)\widehat{B}_r\de r$, and how small $\vB(T)$ is depends on the speeds and weights.
\end{remark}

\section{Scope and Limitations of the Reduction}
\label{sec:limits}

The reduction rests on three properties: The data enter linearly, the noise is Gaussian and independent of the data, and the loss or the cost depends on the path only through a linear statistic.
In this section, we describe the situations in which one of these properties fails, so memory in the noise can change the problem itself and not only its schedule.

\paragraph{Discretization, Parametrization, and Optimization.}
Our results concern continuous-time models, population losses, and exact scores.
In this idealized setting, a fractional model and its Brownian twin share the schedule and, by Corollary~\ref{cor:score-matching}, the regression problem solved in training.
They also generate the same distribution, even when the reverse-time dynamics are started from a Gaussian law at time $T$; see Proposition~\ref{prop:gfdm-initial} for GFDM, and note that the same argument applies to the model of Section~\ref{sec:fbm}, where the statistic at time $T$ is $X_T/c(T)$.
In practice, the two models can still differ.
The integrator of the augmented reverse-time equation and the integrator of~\eqref{eq:twin-reverse} have different discretization errors, and the training dynamics can depend on the input representation.
A fractional model and a conventional, non-twin Brownian model differ in addition in their schedules.
The effective schedule $v$ is determined by the parameters of the noise and is none of the schedules commonly designed for Brownian models, and its noise level at time $T$ can be smaller than the marginal variance suggests (Remark~\ref{rem:gfdm-terminal}).
A comparison of a fractional model with the Brownian model with the effective schedule $v$ is therefore the natural control experiment for any reported effect of memory.

\paragraph{Riemann--Liouville Noise.}
The Markovian approximation of Section~\ref{sec:gfdm} approximates the Riemann--Liouville process $\Gamma(H+\frac12)^{-1}\int_0^t(t-r)^{H-1/2}\de B_r$ and not the fractional Brownian motion of Section~\ref{sec:fbm}.
If this process drives the forward equation~\eqref{eq:forward}, then the observations are still linear-Gaussian.
Writing $G$ for the noise in units of the data, Proposition~\ref{prop:sufficiency} and Theorem~\ref{thm:twin}~\ref{twin:law} hence apply with
\[
    v(t)=\Var\bigl[G^{(1)}_T\given G^{(1)}_u-G^{(1)}_T,\ u\in[t,T]\bigr],\qquad t\in(0,T],
\]
and so does Theorem~\ref{thm:twin}~\ref{twin:sde}, provided that $v$ is positive on $(0,T]$ and absolutely continuous.
The explicit formulas of Theorem~\ref{thm:fbm}, however, rely on the fact that the reversed increments of a fractional Brownian motion are again a fractional Brownian motion with a Volterra representation which generates the same filtration.
Whether the reversed Riemann--Liouville process admits such a representation is left open.

\appendix

\section{Gaussian and Filtering Tools}
\label{sec:auxiliary}

This appendix collects the auxiliary results used in the proofs.
The first one is classical; see, e.g.,~\cite{janson1997gaussian}.

\begin{lemma}[Gaussian Hilbert Spaces]
    \label{lem:gaussian-hilbert}
    Let $\mathcal{G}$ be a centered Gaussian family of real random variables, and let $\mathcal{H}$ be its closed linear span in $\Lp^2(\Prob)$.
    \begin{enumerate}
        \item \label{gh:gaussian}
              The space $\mathcal{H}$ is a centered Gaussian family.
        \item \label{gh:independent}
              If $\mathcal{H}_1,\mathcal{H}_2\subseteq\mathcal{H}$ are orthogonal, then $\sigma(\mathcal{H}_1)$ and $\sigma(\mathcal{H}_2)$ are independent.
        \item \label{gh:projection}
              If $P$ is the orthogonal projection onto a closed subspace $\mathcal{H}_1\subseteq\mathcal{H}$ and $\zeta\in\mathcal{H}$, then $\zeta-P\zeta$ is independent of $\sigma(\mathcal{H}_1)$ and $\E[\zeta\given\sigma(\mathcal{H}_1)]=P\zeta$.
    \end{enumerate}
\end{lemma}

The second result is the innovation argument of filtering theory; see, e.g.,~\cite{liptser2001statistics}.

\begin{lemma}[Projection of Drifts]
    \label{lem:innovation}
    Let $(\Filt_s)_{s\in[0,S]}$ and $(\FiltAlt_s)_{s\in[0,S]}$ be filtrations with $\Filt_s\subseteq\FiltAlt_s$, and let $\zeta$ be an $\Filt$-adapted process with
    \[
        \zeta_s-\zeta_0=\int_0^sb_r\de r+A_s,\qquad s\in[0,S],
    \]
    where $b$ is $\FiltAlt$-progressively measurable with $\E\int_0^S\Abs{b_r}\de r<\infty$ and $A$ is a $\FiltAlt$-martingale.
    If $\hat{b}$ is an $\Filt$-progressively measurable process with $\hat{b}_r=\E[b_r\given\Filt_r]$ almost surely for almost every $r\in[0,S]$, then $\zeta_s-\zeta_0-\int_0^s\hat{b}_r\de r$, $s\in[0,S]$, is an $\Filt$-martingale.
\end{lemma}

\begin{proof}
    The process $\hat{A}_s\defined\zeta_s-\zeta_0-\int_0^s\hat{b}_r\de r=A_s+\int_0^s(b_r-\hat{b}_r)\de r$ is $\Filt$-adapted and, by Jensen's inequality, integrable.
    For $s\leq s'$ and $F\in\Filt_s\subseteq\FiltAlt_s$ we have $\E[(A_{s'}-A_s)\Ind{F}]=0$, and Fubini's theorem gives
    \[
        \E\Bigl[\Ind{F}\int_s^{s'}(b_r-\hat{b}_r)\de r\Bigr]
        =\int_s^{s'}\E\bigl[\Ind{F}\bigl(b_r-\E[b_r\given\Filt_r]\bigr)\bigr]\de r=0,
    \]
    because $F\in\Filt_r$ for $r\geq s$.
    Hence $\E[(\hat{A}_{s'}-\hat{A}_s)\Ind{F}]=0$, and the proof is complete.
\end{proof}

The third result identifies the semimartingale decomposition of a time-changed Brownian motion run backwards in a filtration which may be larger than its own.

\begin{lemma}[Gaussian Bridges]
    \label{lem:bridge}
    Let $\rho\colon[0,T]\to[0,\infty)$ be nonincreasing and absolutely continuous with $\rho>0$ on $[0,T)$ and $\rho(T)=0$.
    Let $\zeta$ be a continuous process in $\R^d$, adapted to a filtration $(\FiltAlt_s)_{s\in[0,T]}$, whose coordinates are independent centered Gaussian processes with covariance $\Cov[\zeta^{(i)}_s,\zeta^{(i)}_{s'}]=\rho(s\vee s')$, and suppose that
    \begin{equation}
        \label{eq:bridge-condition}
        \E\bigl[\zeta_{s'}\given\FiltAlt_s\bigr]=\frac{\rho(s')}{\rho(s)}\,\zeta_s,\qquad 0\leq s\leq s'<T.
    \end{equation}
    Then $\E\int_0^T\Abs{\rho'(r)/\rho(r)}\Abs{\zeta_r}\de r\leq2\sqrt{d\rho(0)}$, the quadratic covariations of $\zeta$ are given by $[\zeta^{(i)},\zeta^{(j)}]_s=\delta_{ij}(\rho(0)-\rho(s))$, and
    \[
        A_s\defined\zeta_s-\zeta_0-\int_0^s\frac{\rho'(r)}{\rho(r)}\,\zeta_r\de r,\qquad s\in[0,T],
    \]
    is a continuous $\FiltAlt$-martingale with $\langle A^{(i)},A^{(j)}\rangle_s=\delta_{ij}(\rho(0)-\rho(s))$.
\end{lemma}

\begin{proof}
    Since $\E\Abs{\zeta_r}\leq\sqrt{d\rho(r)}$ and $\sqrt{\rho}$ is absolutely continuous on $[0,T-\eps]$ for every $\eps>0$, the integral in the first claim is at most $\sqrt{d}\int_0^T(-\rho'(r))\rho(r)^{-1/2}\de r=2\sqrt{d\rho(0)}$.
    For $0\leq s\leq s'<T$, conditional Fubini and~\eqref{eq:bridge-condition} give
    \[
        \E\bigl[A_{s'}-A_s\given\FiltAlt_s\bigr]
        =\frac{\zeta_s}{\rho(s)}\Bigl(\rho(s')-\rho(s)-\int_s^{s'}\rho'(r)\de r\Bigr)=0,
    \]
    and the martingale property extends to $s'=T$, because $A_{s'}\to A_T$ in $\Lp^1$ by the integrability just shown and by $\E\Abs{\zeta_{s'}}^2=d\rho(s')\to0$.
    Since $\rho(s\vee s')=\rho(s)\wedge\rho(s')$, each coordinate of $\zeta$ has the covariance of $(W_{\rho(s)})_{s\in[0,T]}$ for a Brownian motion $W$, so its increments over disjoint intervals are independent with $\E[(\zeta^{(i)}_{s'}-\zeta^{(i)}_s)^2]=\rho(s)-\rho(s')$ for $s\leq s'$.
    For a partition of $[0,s]$, the sum of the squared increments of $\zeta^{(i)}$ thus has mean $\rho(0)-\rho(s)$ and a variance of at most $2(\rho(0)-\rho(s))$ times the largest increment of $\rho$ over the partition, and the sums of products of increments of two different coordinates are centered with a variance bounded in the same way.
    By the uniform continuity of $\rho$, this proves the claim on the quadratic covariations of $\zeta$, and $A$ has the same brackets, since $A-\zeta$ is continuous and of finite variation on $[0,s]$ for every $s<T$.
\end{proof}

The last result of this appendix is used for the reverse-time dynamics of the path in Theorem~\ref{thm:fbm}~\ref{fbm:path-sde}.
It describes the enlargement of a Brownian filtration by a Gaussian functional at time zero and is a special case of the classical results of~\citet{jeulin1980semimartingales}.

\begin{lemma}[Initial Enlargement by a Gaussian Functional]
    \label{lem:enlargement}
    Let $W$ be a Brownian motion in $\R^d$ on $[0,T]$ with augmented natural filtration $\Filt^W$, let $k\in\Lp^2(0,T)$ with $\rho(s)\defined\int_s^Tk(r)^2\de r>0$ for $s<T$, and let $U$ be a random vector independent of $W$.
    Set
    \[
        \Phi\defined\int_0^Tk(r)\de W_r,
        \qquad
        R_s\defined\int_s^Tk(r)\de W_r,
        \qquad
        \FiltAlt_s\defined\sigma(U,\Phi)\vee\Filt^W_s.
    \]
    Then $\int_0^s\Abs{k(r)}\,\E\Abs{R_r}/\rho(r)\de r<\infty$ for $s<T$, and
    \[
        \beta_s\defined W_s-\int_0^s\frac{k(r)}{\rho(r)}\,R_r\de r,\qquad s\in[0,T),
    \]
    is a $\FiltAlt$-Brownian motion.
\end{lemma}

\begin{proof}
    We first show that $\beta$ is well defined.
    For $r\leq s<T$ we have $\E\Abs{R_r}\leq\sqrt{d\rho(r)}$ and $\rho(r)\geq\rho(s)>0$, so the integrability follows from $k\in\Lp^2(0,T)$.
    Since $R_r=\Phi-\int_0^rk\de W$ is continuous and $\FiltAlt$-adapted, so is $\beta$.

    We next reduce the claim to the smaller filtration $\FiltAlt'_s\defined\sigma(\Phi)\vee\Filt^W_s$.
    The random vector $U$ is independent of $(\beta,\Phi,W)$, so an increment of $\beta$ that is independent of $\FiltAlt'_s$ is also independent of $\FiltAlt_s=\sigma(U)\vee\FiltAlt'_s$.

    Fix $s\leq s'<T$ and a coordinate $i$.
    The increment $\beta^{(i)}_{s'}-\beta^{(i)}_s$ is an $\Lp^2$-limit of linear combinations of values of $W$, so it lies, together with $\Phi$ and $(W_q)_{q\leq s}$, in the Gaussian space of $W$.
    It is uncorrelated with $W_q$ for $q\leq s$, because $W_{s'}-W_s$ and $R_r$, $r\geq s$, are.
    It is also uncorrelated with $\Phi^{(j)}$ for $j\neq i$, since different coordinates are independent.
    The key step is the covariance with $\Phi^{(i)}$.
    Since $\Cov[W^{(i)}_{s'}-W^{(i)}_s,\Phi^{(i)}]=\int_s^{s'}k(r)\de r$ and $\Cov[R^{(i)}_r,\Phi^{(i)}]=\rho(r)$, we obtain
    \[
        \Cov\bigl[\beta^{(i)}_{s'}-\beta^{(i)}_s,\Phi^{(i)}\bigr]
        =\int_s^{s'}k(r)\de r-\int_s^{s'}\frac{k(r)}{\rho(r)}\,\rho(r)\de r=0.
    \]
    Each coordinate of the increment is thus orthogonal to $\Phi$ and to $(W_q)_{q\leq s}$, and Lemma~\ref{lem:gaussian-hilbert}~\ref{gh:independent} shows that the increment is independent of $\FiltAlt'_s$.

    Hence $\beta$ is a continuous $\FiltAlt$-martingale on $[0,s']$ for every $s'<T$.
    Since $\beta-W$ is continuous and of finite variation on $[0,s']$, we have $\langle\beta^{(i)},\beta^{(j)}\rangle_s=\langle W^{(i)},W^{(j)}\rangle_s=\delta_{ij}s$, and L\'{e}vy's characterization shows that $\beta$ is a $\FiltAlt$-Brownian motion on $[0,T)$.
\end{proof}

\section{The Molchan--Golosov Kernel}
\label{sec:kernel}

We record the kernel of Lemma~\ref{lem:molchan-golosov}, following~\cite[Section~5.1.3]{nualart2006malliavin}.
For $H>\frac12$ and $0<r<t$, it is given by
\[
    K_H(t,r)=c_Hr^{\frac12-H}\int_r^t(u-r)^{H-\frac32}u^{H-\frac12}\de u,
    \qquad
    c_H\defined\Bigl(\frac{H(2H-1)}{\mathrm{B}(2-2H,H-\frac12)}\Bigr)^{1/2},
\]
and for $H<\frac12$ and $0<r<t$ by
\begin{align}
    K_H(t,r) & =c_H\Bigl[\Bigl(\frac{t}{r}\Bigr)^{H-\frac12}(t-r)^{H-\frac12}-\bigl(H-\tfrac12\bigr)r^{\frac12-H}\int_r^tu^{H-\frac32}(u-r)^{H-\frac12}\de u\Bigr], \\
    c_H      & \defined\Bigl(\frac{2H}{(1-2H)\,\mathrm{B}(1-2H,H+\frac12)}\Bigr)^{1/2},
\end{align}
where $\mathrm{B}$ is the Beta function.
For $H=\frac12$, the kernel is $K_H\equiv1$.
We set $c'_H\defined c_H/(H-\frac12)$ for $H>\frac12$, $c'_H\defined c_H$ for $H<\frac12$, and $c'_H\defined1$ for $H=\frac12$.

\begin{lemma}[Properties of the Kernel]
    \label{lem:kernel}
    Let $H\in(0,1)$.
    Then the following hold.
    \begin{enumerate}
        \item \label{kernel:positive}
              The kernel is positive, that is, $K_H(t,r)>0$ for $0<r<t$.
        \item \label{kernel:bound}
              For every $\delta\in(0,T)$ there is a constant $C$ such that $K_H(u,r)\leq C(u-r)^{H-\frac12}$ for $\delta\leq r<u\leq T$.
        \item \label{kernel:asymptotics}
              For every $t>0$, it holds that $K_H(t,r)/(t-r)^{H-\frac12}\to c'_H$ as $r\uparrow t$.
        \item \label{kernel:kappa}
              The constant $c'_H$ satisfies $(c'_H)^2/(2H)=\Gamma(\frac32-H)/(\Gamma(H+\frac12)\Gamma(2-2H))$.
    \end{enumerate}
\end{lemma}

\begin{proof}
    For $H=\frac12$, all statements are trivial, so let $H\neq\frac12$.

    We first treat the case $H>\frac12$.
    Let $\delta\leq r<u\leq T$.
    In the integral defining $K_H(u,r)$, the variable of integration $q$ satisfies $r\leq q\leq u$, hence $r^{H-1/2}\leq q^{H-1/2}\leq u^{H-1/2}$.
    Moreover, we have $\int_r^u(q-r)^{H-3/2}\de q=(u-r)^{H-1/2}/(H-\frac12)$.
    This yields the two-sided bound
    \[
        c'_H(u-r)^{H-\frac12}
        \leq K_H(u,r)
        \leq c'_H\Bigl(\frac{u}{r}\Bigr)^{H-\frac12}(u-r)^{H-\frac12}.
    \]
    The lower bound gives~\ref{kernel:positive}, and the upper bound gives~\ref{kernel:bound} with $C=c'_H(T/\delta)^{H-1/2}$.
    Since $(u/r)^{H-1/2}\to1$ as $r\uparrow u$, the two bounds together give~\ref{kernel:asymptotics}.

    We now treat the case $H<\frac12$.
    Both terms in the bracket are positive, since $(t/r)^{H-1/2}>0$ and $-(H-\frac12)>0$, which gives~\ref{kernel:positive}.
    Let $\delta\leq r<u\leq T$.
    The first term in the bracket is at most $(u-r)^{H-1/2}$, because $(u/r)^{H-1/2}\leq1$.
    In the second term, we use $q^{H-3/2}\leq r^{H-3/2}$ for $q\geq r$, which shows that this term is at most
    \[
        \bigl(\tfrac12-H\bigr)\,r^{\frac12-H}r^{H-\frac32}\,\frac{(u-r)^{H+\frac12}}{H+\frac12}
        \leq\delta^{-1}(u-r)^{H+\frac12}
        \leq\delta^{-1}T(u-r)^{H-\frac12}.
    \]
    This gives~\ref{kernel:bound} with $C=c_H(1+T/\delta)$.
    Moreover, as $r\uparrow u$, the first term divided by $(u-r)^{H-1/2}$ equals $(u/r)^{H-1/2}\to1$, while the second term divided by $(u-r)^{H-1/2}$ is at most $\delta^{-1}(u-r)\to0$.
    Since $c'_H=c_H$ in this case, this gives~\ref{kernel:asymptotics}.

    Finally, \ref{kernel:kappa} follows by a direct calculation from the definitions of $c_H$ and $c'_H$.
    It uses the relation $\mathrm{B}(a,b)=\Gamma(a)\Gamma(b)/\Gamma(a+b)$ and the functional equation $x\Gamma(x)=\Gamma(x+1)$, the latter with $x=H-\frac12$ for $H>\frac12$ and with $x=1-2H$ for $H<\frac12$.
\end{proof}

\printbibliography

\end{document}